\documentclass[runningheads,5pt]{article}
\usepackage{amssymb}
\usepackage{eucal}
\usepackage[dvips]{graphicx}
\usepackage{hyperref}
\usepackage{multirow}
\usepackage{fancyhdr}
\renewcommand{\baselinestretch}{1.2}
\begin{document}

\title{Reaction Diffusion Equations with Nonlinear Boundary Conditions in Narrow Domains}

\author{Mark Freidlin\footnotemark[2] \hspace{0.1cm} and Konstantinos
Spiliopoulos\footnotemark[3]\\
February 2008}

\date{}
\maketitle
$^{\dag}$\textit{Department of Mathematics, University of
Maryland, College Park, 20742, }

$ $ \textit{Maryland, USA, E-mail: mif@math.umd.edu}

$^{\ddag}$\textit{Department of Mathematics, University of Maryland,
College Park, 20742,}

$ $ \textit{Maryland, USA, kspiliop@math.umd.edu}

\begin{abstract}
Second initial boundary problem in narrow domains of width $\epsilon\ll 1$ for linear second order differential
equations with nonlinear boundary conditions is considered in this paper. Using probabilistic methods we show that the solution of such a problem converges as $\epsilon \downarrow 0$ to the solution of a standard reaction-diffusion equation in a domain of reduced dimension. This reduction allows to obtain some results concerning wave front
propagation in narrow domains. In particular, we describe conditions leading to jumps of the wave front.

Keywords: Reaction Diffusion equations, Narrow Domains, Wave Front
Propagation, Instantaneous Reflection
\end{abstract}

\vspace{2cm}

\framebox{\parbox[c]{11cm}{This is an electronic reprint of the
original article published by the
\href{http://www.iospress.nl/loadtop/load.php?isbn=09217134}{
Asymptotic Analysis Volume 59, Number 3-4 / 2008, pp.227-249}. This
reprint differs from the original in pagination, typographic detail
and typo corrections.}}

\newpage

\section{Introduction}

 For each $x\in \mathbb{R}^{n}$, let
$D_{x}$ be a bounded domain in $\mathbb{R}^{m}$ with a smooth
boundary $\partial D_{x}$. Assume, for brevity, that $D_{x}$ is
homeomorphic to a ball in $\mathbb{R}^{m}$ and contains
$0\in\mathbb{R}^{m}$. Consider the domain $D=\{(x,y):x\in
\mathbb{R}^{n},y\in D_{x}\}\subset \mathbb{R}^{n+m}$. Assume that
the boundary $\partial D$ of $D$ is smooth enough and denote by
$\gamma(x,y)$ the inward unit normal to $\partial D$.  Assume that
$\gamma(x,y)$ is not parallel to the subspace $\mathbb{R}^{n}\subset
\mathbb{R}^{n+m}$ for any $(x,y)\in \partial D$.

Denote by $D^{\epsilon}, 0<\epsilon<<1,$ the domain in  $\mathbb{R}^{n+m}$ obtained from $D$ by contraction: $D^{\epsilon}=\{(x,y): x\in \mathbb{R}^{n}, y\epsilon^{-1}\in D_{x}\}$. If $n=1$, $D^{\epsilon}$ is a narrow tube (or a strip for $m=1$) for $0<\epsilon<<1$. If $n>1$, then $D^{\epsilon}$ is a thin layer.

Consider the problem:
\begin{eqnarray}
u^{\epsilon}_{t}&=&\frac{1}{2}\triangle u^{\epsilon}, \hspace{2cm} \textrm{in}\hspace{0.1cm} (0,T) \times D^{\epsilon}\label{InitialEquationWith_e}\\
u^{\epsilon}(0,x,y)&=&f(x), \hspace{2.2cm} \textrm{on}\hspace{0.1cm} \{0\} \times D^{\epsilon}\nonumber\\
\frac{\partial u^{\epsilon}}{\partial \gamma^{\epsilon}}&=&-\epsilon
c(x,y,u^{\epsilon})u^{\epsilon}, \hspace{0.7cm} \textrm{on}
\hspace{0.1cm}(0,T) \times \partial D^{\epsilon},\nonumber
\end{eqnarray}
where $\gamma^{\epsilon}$ is the inward unit normal to $\partial
D^{\epsilon}$. The functions $f$ and $c$ are
sufficiently regular and bounded; $f$ is assumed to be nonnegative. Our goal in this paper is
to study the behavior of solution of problem (\ref{InitialEquationWith_e}) as $\epsilon\downarrow 0$. Using probabilistic methods, we will prove that
$u^{\epsilon}(t,x,y)$ converges as $\epsilon\downarrow 0$ to the solution $u(t,x)$ of the problem:
\begin{eqnarray}
u_{t}&=&\frac{1}{2}\triangle_{x} u +\frac{1}{2} \nabla (\log V(x))\nabla_{x} u + \frac{1}{2}\frac{S(x)}{V(x)}c(x,0,u)u, \hspace{0.3cm} \textrm{in}\hspace{0.1cm} (0,T) \times \mathbb{R}^{n}\nonumber \\
u(0,x)&=&f(x), \hspace{1.1cm} \textrm{on}\hspace{0.1cm} \{0\} \times \mathbb{R}^{n}.\label{InitialEquationWithout_e}
\end{eqnarray}
Here $V(x)$ is the volume of $D_{x}$ and $S(x)$ is the surface area of $\partial D_{x}$. One can expect that, under certain assumptions on the nonlinear term  $c(x,y,u)u$ in (\ref{InitialEquationWith_e}), the solution $u^{\epsilon}(t,x,y)$ can be approximated by a running-wave-type solution. Corresponding results on the standard reaction diffusion equation (\ref{InitialEquationWithout_e}) (see chapter $6$ and $7$ in \cite{K1}) allow to describe the asymptotic wavefront motion for (\ref{InitialEquationWith_e}). We will see how the motion of the interface (wavefront) depends on the behavior of the cross-sections $D_{x}$ of the domain $D$. In particular, using the results of \cite{K1} (chapter 6) we will see that in the case of the nonlinear term of K-P-P type the wavefront can have jumps.

Consider the Wiener process $(X^{\epsilon}_{t},Y^{\epsilon}_{t})$ in $D^{\epsilon}$ with instantaneous normal reflection on the boundary of $D^{\epsilon}$. Its trajectories can be described by the stochastic differential equations:
\begin{eqnarray}
X_{t}^{\epsilon}&=& x+ W_{t}^{1}+\int_{0}^{t}\gamma_{1}^{\epsilon}(X_{s}^{\epsilon},Y_{s}^{\epsilon})dL_{s}^{\epsilon}\nonumber\\
Y_{t}^{\epsilon}&=& y + W_{t}^{2}+\int_{0}^{t}\gamma_{2}^{\epsilon}(X_{s}^{\epsilon},Y_{s}^{\epsilon})dL_{s}^{\epsilon}.\label{StochasticProcessWithReflection1}\renewcommand{\baselinestretch}{1.3}
\end{eqnarray}
Here  $W_{t}^{1}$ and $W_{t}^{2}$ are independent Wiener process in $\mathbb{R}^{n}$ and $\mathbb{R}^{m}$ respectively and $(x,y)$ is a point inside $D^{\epsilon}$. Moreover $\gamma_{1}^{\epsilon}$ and $\gamma_{2}^{\epsilon}$ are projections of  the unit inward normal vector to $\partial D^{\epsilon}$ on $\mathbb{R}^{n}$ and $\mathbb{R}^{m}$ respectively. It is easy to see that $\lim_{\epsilon\downarrow 0}|\epsilon^{-1}\gamma_{1}^{\epsilon}|=\frac{\gamma_{1}^{1}}{|\gamma_{2}^{1}|}$ and $\lim_{\epsilon\downarrow 0}|\gamma_{2}^{\epsilon}|=1$, where $|\cdot|$ denotes Euclidean length. Furthermore $L^{\epsilon}_{t}$ is the local time for the process $(X^{\epsilon}_{t},Y^{\epsilon}_{t})$ on $\partial D^{\epsilon}$, i.e. it is a continuous, non-decreasing process that increases only when $(X^{\epsilon}_{t},Y^{\epsilon}_{t}) \in \partial D^{\epsilon}$ such that the Lebesque measure $\Lambda\{t>0:(X^{\epsilon}_{t},Y^{\epsilon}_{t}) \in \partial D^{\epsilon}\}=0$ (see for instance \cite{KS}).

If $(X^{\epsilon}_{t},Y^{\epsilon}_{t})$ is defined by (\ref{StochasticProcessWithReflection1}), then as it can be derived from Theorem 2.5.1 in \cite{K1}, $u^{\epsilon}(t,x,y)$ satisfies the following integral equation in the functional space:
\begin{equation}
u^{\epsilon}(t,x,y)=E_{x,y}f(X^{\epsilon}_{t})\exp[\int_{0}^{t}\epsilon
c(X^{\epsilon}_{s},Y^{\epsilon}_{s},u^{\epsilon}(t-s,X^{\epsilon}_{s},Y^{\epsilon}_{s}))dL^{\epsilon}_{s}],\label{FeynmanKacFormulaBeforeLimit1}
\end{equation}
where $E_{x,y}$ denotes expectation and the subscript $(x,y)$ denotes the initial point of $(X^{\epsilon}_{s},Y^{\epsilon}_{s})$. Equation (\ref{FeynmanKacFormulaBeforeLimit1}) has a unique solution if, say, $c(x,y,u)$ has a bounded derivative in $u$.

Let $X_{t}$ be the solution of the stochastic differential equation
\begin{equation}
X_{t}= x + W_{t}^{1} + \int_{0}^{t}\frac{1}{2}\nabla (\log V(X_{s}))ds. \label{LimitingStochasticProcess1}
\end{equation}

Then the solution $u(t,x)$ of equation (\ref{InitialEquationWithout_e}) satisfies the equality:
\begin{equation}
u(t,x)=E_{x}f(X_{t})\exp[\int_{0}^{t}\frac{1}{2}\frac{S(X_{s})}{V(X_{s})}c(X_{s},0,u(t-s,X_{s}))ds].\label{FeynmanKacFormulaAfterLimit1}
\end{equation}

We prove that the component $X^{\epsilon}_{t}$ of the process
$(X^{\epsilon}_{t},Y^{\epsilon}_{t})$ converges in a certain sense
to $X_{t}$. This together with uniform in $0<\epsilon<1$ bounds for
$u^{\epsilon}(t,x,y)$ and its derivatives allow to prove that the
solution of (\ref{FeynmanKacFormulaBeforeLimit1}) converges to the
solution of (\ref{FeynmanKacFormulaAfterLimit1}) as
$\epsilon\downarrow 0$ uniformly on each compact subset of
$[0,\infty)\times \mathbb{R}^{n+m}$.

In the next section we consider averaging of integrals in local
time. This result allows in section 3 to prove convergence of the
integral in the right side of the first of equations in
(\ref{StochasticProcessWithReflection1}) to the integral term in
(\ref{LimitingStochasticProcess1}) and convergence of exponents in
(\ref{FeynmanKacFormulaBeforeLimit1}) and
(\ref{FeynmanKacFormulaAfterLimit1}). Together with a-priori  bounds
obtained in section 3, this implies convergence of
$u^{\epsilon}(t,x,y)$ to $u(t,x)$. Some results concerning wavefront
propagation are presented in section 4.

\section{Averaging of Integrals in Local Time}

Let $H(x,y)$ be a smooth and bounded function. We want to consider
the limiting behavior as $\epsilon\downarrow 0$ of expressions like
$\int_{0}^{t}\epsilon H(X^{\epsilon}_{s},Y^{\epsilon}_{s}/\epsilon)d
L^{\epsilon}_{s}$ (see Lemma 2.1 below). We will assume that the
unit inward normal $\gamma(x,y)$ to $\partial D$ and the function
$H(x,y)$ are both three times differentiable in $x$ and $y$.

\vspace{0.2cm}

\textbf{Lemma 2.1.} Define $Q(x)=\frac{1}{V(x)}\int_{\partial D_{x}}H(x,y)dS_{x}$, where $dS_{x}$ is the surface element on $\partial D_{x}$. Then for every $T>0$ and small enough $\epsilon$,
there exists a constant $K$ independent of $\epsilon$ such that:
\begin{enumerate}
\item{$\sup_{0\leq t \leq T}E|\int_{0}^{t}\frac{1}{2}Q(X_{s}^{\epsilon})ds-\int_{0}^{t}\epsilon H(X_{s}^{\epsilon},Y_{s}^{\epsilon}/\epsilon)|\gamma_{2}^{\epsilon}(X_{s}^{\epsilon},Y_{s}^{\epsilon})|dL^{\epsilon}_{s} |^{2}\leq K \epsilon^{2}$.}
\item{For every $\delta>0$ we have

 $P\{\sup_{0\leq t \leq T}|\int_{0}^{t}\frac{1}{2}Q(X_{s}^{\epsilon})ds-\int_{0}^{t}\epsilon H(X_{s}^{\epsilon},Y_{s}^{\epsilon}/\epsilon)|\gamma_{2}^{\epsilon}(X_{s}^{\epsilon},Y_{s}^{\epsilon})|dL^{\epsilon}_{s} |>\delta\}\leq K\frac{\epsilon^{2}}{\delta^{2}}$.}
\end{enumerate}

The proof of lemma 2.1 relies on the following lemma, which we prove first.

\vspace{0.2cm}

\textbf{Lemma 2.2.} For every $T>0$ and small enough $\epsilon$,
there exists a constant $K_{1}$ independent of $\epsilon$ such that:
\begin{displaymath}
E |\epsilon L_{T}^{\epsilon}|^{2}\leq K_{1}
\end{displaymath}\label{LemmaRDE22}

\textbf{Proof.}

Consider the auxiliary problem
\begin{eqnarray}
\triangle_{y} v(x,y)&=&Q(x), \hspace{0.2cm} y\in D_{x}\subset \mathbb{R}^{m}\nonumber\\
\frac{\partial_{y} v(x,y)}{\partial n(x,y)}&=&-1, \hspace{0.2cm} y\in \partial D_{x}, \label{SpecificNeymmanProblem}
\end{eqnarray}
where $n(x,y)=\frac{\gamma_{2}^{1}(x,y)}{|\gamma_{2}^{1}(x,y)|}$ and $ x \in \mathbb{R}^{n}$ is a parameter. Let
\begin{equation}
Q(x)=\frac{S(x)}{V(x)},\label{SpecificSolvabilityCondition}
\end{equation}
where $S(x)$ is the surface area of $D_{x}$ and $V(x)$ is the volume of $D_{x}$. As it can be derived from \cite{BM}, a smooth in $x$ and $y$ solution $v(x,y)$ of problem (\ref{SpecificNeymmanProblem}) exists and is bounded together with its first and second derivatives. So we can apply It\^{o} formula to the function $\epsilon v(x,y/\epsilon)$, and get:
\begin{eqnarray}
\epsilon^{2}v(X_{t}^{\epsilon},Y_{t}^{\epsilon}/\epsilon)&=&\epsilon^{2}v(x,y/\epsilon)+
\int_{0}^{t}\epsilon^{2}\frac{1}{2}\triangle_{x}v(X_{s}^{\epsilon},Y_{s}^{\epsilon}/\epsilon)ds+\int_{0}^{t}\frac{1}{2}\triangle_{y}v(X_{s}^{\epsilon},Y_{s}^{\epsilon}/\epsilon)ds\nonumber\\
&+&\int_{0}^{t}\epsilon^{2}(\nabla_{x}v(X_{s}^{\epsilon},Y_{s}^{\epsilon}/\epsilon),dW_{s}^{1})+\int_{0}^{t}\epsilon(\nabla_{y}v(X_{s}^{\epsilon},Y_{s}^{\epsilon}/\epsilon),dW_{s}^{2})\nonumber\\
&+&\int_{0}^{t}\epsilon^{2}(\nabla_{x}v(X_{s}^{\epsilon},Y_{s}^{\epsilon}/\epsilon),\gamma^{\epsilon}_{1}(X_{s}^{\epsilon},Y_{s}^{\epsilon}))dL^{\epsilon}_{s}\nonumber\\
&+&\int_{0}^{t}\epsilon(\nabla_{y}v(X_{s}^{\epsilon},Y_{s}^{\epsilon}/\epsilon),\gamma^{\epsilon}_{2}(X_{s}^{\epsilon},Y_{s}^{\epsilon}))dL^{\epsilon}_{s}\label{ItoFormulaForLocalTime}
\end{eqnarray}
Recalling now that $\lim_{\epsilon\downarrow 0}|\epsilon^{-1}\gamma_{1}^{\epsilon}|=\frac{\gamma_{1}^{1}}{|\gamma_{2}^{1}|}$ and $\lim_{\epsilon\downarrow 0}|\gamma_{2}^{\epsilon}|=1$ and that $v$ satisfies (\ref{SpecificNeymmanProblem}) one easily concludes that there is an
$\epsilon_{0}=\epsilon_{0}(\| |\nabla_{x} v| \|,\gamma_{1}^{1})>0$ such that for every $\epsilon < \epsilon_{0}$:
\begin{eqnarray}
E|\epsilon L_{T}^{\epsilon}|^{2}&\leq& C[\epsilon^{4}(2\|v^{2}\|+\|\frac{1}{2}\triangle_{x}v\|^{2}T^{2}+\||\nabla_{x}v|^{2}\|T)+\nonumber\\
&+&\epsilon^{2}\||\nabla_{y}v|^{2}\|T+\|\frac{1}{2}Q\|^{2}T],\nonumber
\end{eqnarray}
where for any function $g$, $\| g \|=\sup_{z} |g(z)|$. Here, we also
used the fact that the local time is increasing function of $t$.
This proves the statement of the lemma.
%
\begin{flushright}
$\square$
\end{flushright}

\textbf{Proof of Lemma 2.1.}
We consider the auxiliary problem
\begin{eqnarray}
\triangle_{y} v(x,y)&=&Q(x), \hspace{0.2cm} y\in D_{x}\subset \mathbb{R}^{m}\nonumber\\
\frac{\partial_{y} v(x,y)}{\partial n(x,y)}&=&-H(x,y), \hspace{0.2cm} y\in \partial D_{x}, \label{GeneralNeymmanProblem1}
\end{eqnarray}
where $n(x,y)=\frac{\gamma_{2}^{1}(x,y)}{|\gamma_{2}^{1}(x,y)|}$ and $ x \in \mathbb{R}^{n}$ is a parameter.

The necessary and sufficient condition for the existence of a solution for (\ref{GeneralNeymmanProblem1}) is that
\begin{equation}
Q(x)=\frac{1}{V(x)}\int_{\partial D_{x}}H(x,y)dS_{x},\label{GeneralSolvabilityCondition}
\end{equation}
where $dS_{x}$ is the surface element on $\partial D_{x}$ and $V(x)=\textrm{volume}(D_{x})$.

Applying It\^o formula to the function $\epsilon v(x,y/\epsilon)$ and using the bounds obtained in Lemma 2.2 we get the following inequalities:
\begin{eqnarray}
& &\sup_{0\leq t \leq
T}E|\int_{0}^{t}\frac{1}{2}Q(X_{s}^{\epsilon})ds-\int_{0}^{t}
\epsilon
H(X_{s}^{\epsilon},Y_{s}^{\epsilon}/\epsilon)|\gamma_{2}^{\epsilon}(X_{s}^{\epsilon},Y_{s}^{\epsilon})|dL^{\epsilon}_{s}
|^{2}\leq \nonumber\\
&\leq&\epsilon^{4}C(2\|v^{2}\|+\|\frac{1}{2}\triangle_{x}v\|^{2}T^{2}+\||\nabla_{x}v|^{2}\|T+\||\nabla_{x}v|\|^{2}K_{1})+\epsilon^{2}C\||\nabla_{y}v|^{2}\|T,\nonumber
\end{eqnarray}
which proves statement (i) of the lemma.

For part (ii) one makes use of the Doob maximal inequalities (see
\cite{KS}, page 14):
\begin{eqnarray}
E[\sup_{0\leq t \leq T}|\int_{0}^{t}(\nabla_{x}v(X_{s}^{\epsilon},Y_{s}^{\epsilon}/\epsilon),dW_{s}^{1})|]^{2}&\leq& 4 \parallel |\nabla_{x}v|^{2}\parallel T\nonumber\\
E[\sup_{0\leq t \leq
T}|\int_{0}^{t}(\nabla_{y}v(X_{s}^{\epsilon},Y_{s}^{\epsilon}/\epsilon),dW_{s}^{2})|]^{2}&\leq&
4 \parallel |\nabla_{y}v|^{2}\parallel T\nonumber
\end{eqnarray}
Then, following the procedure that proved part $(i)$ we get that
there is an $\epsilon_{0}>0$ such that for every $\epsilon <
\epsilon_{0}$:
\begin{eqnarray}
& &E[\sup_{0\leq t \leq
T}|\int_{0}^{t}\frac{1}{2}Q(X_{s}^{\epsilon})ds-\int_{0}^{t}
\epsilon
H(X_{s}^{\epsilon},Y_{s}^{\epsilon}/\epsilon)|\gamma_{2}^{\epsilon}(X_{s}^{\epsilon},Y_{s}^{\epsilon})|dL^{\epsilon}_{s}
|]^{2}\leq \nonumber\\
&\leq&\epsilon^{4}C(4\|v^{2}\|+\|\frac{1}{2}\triangle_{x}v\|^{2}T^{2}+4\||\nabla_{x}v|^{2}\|T+\||\nabla_{x}v|\|^{2}K_{1})+\epsilon^{2}C\||\nabla_{y}v|^{2}\|T,\nonumber
\end{eqnarray}
which together with Chebyshev inequality proves statement (ii) of
the lemma.
\begin{flushright}
$\square$
\end{flushright}

\section{Limit of $u^{\epsilon}$.}

In this section we consider the limit as $\epsilon\rightarrow 0$ of
the solution $u^{\epsilon}$ to problem
(\ref{InitialEquationWith_e}).
The result is given in Theorem 3.4. The proof  will proceed as
follows. First (in Proposition 3.2.) we write down an integral equation in the space of trajectories for the solution of (\ref{InitialEquationWith_e}). Then
in Lemma 3.3 we consider the mean square limit as
$\epsilon\rightarrow 0$ of the underlying stochastic process with instantaneous normal
reflection on the boundary of $D^{\epsilon}$ (see (\ref{StochasticProcessWithReflection1})). Lastly an important ingredient to the proof are the
a-priori bounds for $u^{\epsilon}$ and its derivatives. These
a-priori bounds are independent of $\epsilon$, their derivation is
standard and are given for completeness in Proposition 3.7.

We assume that the initial function $f(x)$ of problem
(\ref{InitialEquationWith_e}) is bounded, non-negative and can have finite number of simple discontinuities. The function
$c(x,y,u)$ is assumed to be uniformly bounded in all arguments, continuous in x,y, Lipschitz continuous
in u and that there exist constants $M,N>0$ such that $c(\cdot,\cdot,u)<-M$ for $u>N$.

\vspace{0.2cm}

In addition we assume that the boundary of $D^{1}$
satisfies $\partial D^{1} \in
\mathcal{C}^{3+a}(\mathbb{R}^{m})$, where $a \in (0,1)$.

\vspace{0.2cm}

\textbf{Remark 3.1.} For the existence of a classical solution to (\ref{InitialEquationWith_e}) one actually needs only to assume $\partial D^{1} \in \mathcal{C}^{2+a}(\mathbb{R}^{m})$. The assumption $\partial D^{1} \in \mathcal{C}^{3+a}(\mathbb{R}^{m})$ is being done solely for the purpose of Lemma 3.2 and Theorem 3.3.

\vspace{0.2cm}

Let
$(X^{\epsilon},Y^{\epsilon},L^{\epsilon})$ in $\mathbb{R}^{n}\times
\mathbb{R}^{m}\times \mathbb{R}^{1}_{+}$  satisfy
(\ref{StochasticProcessWithReflection1}). Then we have:

\vspace{0.2cm}

\textbf{Proposition 3.2.} Problem (\ref{InitialEquationWith_e}) has a
unique classical solution in $[0,T) \times D^{\epsilon}$ which satisfies:
\begin{equation}
u^{\epsilon}(t,x,y)=E_{x,y}f(X^{\epsilon}_{t})\exp[\int_{0}^{t}\epsilon
c(X^{\epsilon}_{s},Y^{\epsilon}_{s},u^{\epsilon}(t-s,X^{\epsilon}_{s},Y^{\epsilon}_{s}))dL^{\epsilon}_{s}].\label{FeynmanKacFormulaBeforeLimit}
\end{equation}
\textbf{Proof.}  Under our assumptions, the uniqueness and existence of a classical solution to (\ref{InitialEquationWith_e}) follows from Theorem 7.5.13 of \cite{AF1}). The equation (\ref{FeynmanKacFormulaBeforeLimit}) follows from Theorem 2.5.1 of \cite{K1}.
\begin{flushright}
$\square$
\end{flushright}

In order now to consider the  limit as $\epsilon\rightarrow 0$ of (\ref{FeynmanKacFormulaBeforeLimit}), we need first to examine
the asymptotic behavior of $X^{\epsilon}_{t}$ as $\epsilon\rightarrow 0$.

We will prove that $X^{\epsilon}_{t}$ converges as $\epsilon\downarrow 0$ to $X_{t}$, where $X_{t}$ is the solution to
\begin{equation}
X_{t}= x + W_{t}^{1} + \int_{0}^{t}\frac{1}{2}\nabla (\log V(X_{s}))ds, \label{LimitingStochasticProcess}
\end{equation}
where $V(x)=\textrm{volume}(D_{x})$. Hence, we see that as $\epsilon\downarrow 0$, the effect of the reflection on the boundary is an extra drift term. A sketch of the proof for the above result is given in chapter 7 of \cite{F2}. More details are given here.

\vspace{0.2cm}

\textbf{Lemma 3.3.} For any $T>0$ we have
\begin{equation}
\sup_{0\leq t \leq T}E_{x}|X_{t}^{\epsilon}-X_{t}|^{2}\rightarrow 0 \hspace{0.2cm} \textrm{as} \hspace{0.2cm} \epsilon \rightarrow 0.
\end{equation}

\textbf{Proof.}  It is not difficult to see that $\gamma^{\epsilon}_{1}(x,y)=\epsilon\frac{\gamma_{1}^{1}(x,y)}{|\gamma_{2}^{1}(x,y)|}|\gamma_{2}^{\epsilon}(x,y)|$ and that
\begin{equation}
\int_{\partial D_{x}}\frac{\gamma_{1}^{1}(x,y)}{|\gamma_{2}^{1}(x,y)|}dS_{x}=\nabla V(x).\label{DriftTermAndVolume}
\end{equation}
Then Lemma 2.1 with $H(x,y)=\frac{\gamma_{1}^{1}(x,y)}{|\gamma_{2}^{1}(x,y)|}$ and $Q(x)=\nabla \log V(x)$ implies that for small enough $\epsilon$ there exists a constant $K$ independent of $\epsilon$ such that
\begin{equation}
\sup_{0\leq t \leq T}E|\int_{0}^{t}\frac{1}{2}\nabla \log (V(X^{\epsilon}_{s}))ds-\int_{0}^{t}\gamma_{1}^{\epsilon}(X_{s}^{\epsilon},Y_{s}^{\epsilon})dL^{\epsilon}_{s} |^{2}\leq \epsilon^{2}K. \label{LimitLocalTimeIntegral2}
\end{equation}

Now we write
\begin{eqnarray}
X_{t}^{\epsilon}-X_{t}&=& \int_{0}^{t}\gamma_{1}^{\epsilon}(X_{s}^{\epsilon},Y_{s}^{\epsilon})dL^{\epsilon}_{s}-\int_{0}^{t}\frac{1}{2}\nabla \log (V(X_{s}))ds\nonumber\\
&=& [\int_{0}^{t}\gamma_{1}^{\epsilon}(X_{s}^{\epsilon},Y_{s}^{\epsilon})dL^{\epsilon}_{s}-\int_{0}^{t}\frac{1}{2}\nabla \log (V(X^{\epsilon}_{s}))ds]\nonumber\\
&+&[\int_{0}^{t}\frac{1}{2}\nabla \log (V(X^{\epsilon}_{s}))ds-\int_{0}^{t}\frac{1}{2}\nabla \log (V(X_{s}))ds]\label{BeforeGronwallInequalityForStochasticProcess}
\end{eqnarray}

Then Gronwall Lemma
and (\ref{LimitLocalTimeIntegral2}) give:
\begin{equation}
\sup_{0\leq t \leq T}E_{x}|X_{t}^{\epsilon}-X_{t}|^{2}\rightarrow 0 \hspace{0.2cm} \textrm{as} \hspace{0.2cm} \epsilon \rightarrow 0,
\end{equation}
which is the statement of the lemma.
\begin{flushright}
$\square$
\end{flushright}

Consider now the solution, $u$, to the equation
\begin{equation}
u(t,x)=E_{x}f(X_{t})\exp[\int_{0}^{t}
\overline{c}(X_{s},u(t-s,X_{s}))ds],\label{FeynmanKacFormulaAfterLimit}
\end{equation}
where
\begin{equation}
\overline{c}(x,u(t,x))=\frac{1}{2}\frac{S(x)}{V(x)}c(x,0,u(t,x)).\label{NewNonLinearity}
\end{equation}
For notational
convenience we will also denote $\overline{c}(t,x)=\overline{c}(x,u(t,x))$.

Since $\overline{c}(x,u)$ is Lipschitz continuous in $u$, the solution of (\ref{FeynmanKacFormulaAfterLimit}) exists and is unique. Our assumptions on the functions $f,c$ and the boundary $\partial D_{x}$, imply that the solution $u$ to (\ref{FeynmanKacFormulaAfterLimit}) is actually the classical solution of the following parabolic problem:
\begin{eqnarray}
& &u_{t}=\frac{1}{2}\triangle_{x} u +\frac{1}{2} \nabla (\log V(x))\cdot \nabla_{x} u + \frac{1}{2}\frac{S(x)}{V(x)}c(x,0,u(t,x))u, \hspace{0.1cm} \textrm{in}\hspace{0.1cm} (0,T) \times \mathbb{R}^{n}\nonumber\\
& &u(0,x)=f(x), \hspace{0.1cm} \textrm{on}\hspace{0.1cm} \{0\} \times \mathbb{R}^{n}.\label{FinalEquationWithout_e}
\end{eqnarray}

\vspace{0.2cm}

\textbf{Theorem 3.4.} Under our assumptions, we have that
$$u^{\epsilon}(t,x,y)\rightarrow u(t,x)\textrm{ as }\epsilon\rightarrow 0 \textrm{, uniformly in any compact sunset of } \mathbb{R}_{+}\times \mathbb{R}^{n}\times \mathbb{R}^{m},$$
where $u^{\epsilon}(t,x,y)$, $u(t,x)$ are the solutions to (\ref{InitialEquationWith_e})
and (\ref{FinalEquationWithout_e}) respectively.

\textbf{Proof. } By Proposition 3.6 and the well known theorem of Ascoli-Arzela we get that there exists a subsequence of $\{u^{\epsilon}\}$ (which for convenience we will denote again by
$\{u^{\epsilon}\}$) and a function $u$, such that:
$$u^{\epsilon}\rightarrow u\textrm{ as }\epsilon\rightarrow 0 \textrm{, uniformly in compacts}.$$
We will prove that $u$ actually satisfies  (\ref{FeynmanKacFormulaAfterLimit}) which then implies that $u$ satisfies (\ref{FinalEquationWithout_e}). Fix $t$ and $x$ and consider the solution $v(y)=v^{\epsilon,t,x}(y)$ to the elliptic boundary
value problem:
\begin{eqnarray}
\triangle_{y} v(y)&=&\overline{c}^{\epsilon}(t,x), \hspace{0.2cm} y\in D_{x}\subset \mathbb{R}^{m}\nonumber\\
\frac{\partial_{y} v(y)}{\partial n(x,y)}&=&
-\frac{1}{|\gamma_{2}^{\epsilon}(x,\epsilon y)|}c(x,\epsilon y,
u^{\epsilon}(t,x,\epsilon y)), \hspace{0.2cm} y\in \partial D_{x}.
\label{NeymmanProblemForParabolicProblemAndu}
\end{eqnarray}
Problem
(\ref{NeymmanProblemForParabolicProblemAndu}) is solvable if
\begin{equation}
\overline{c}^{\epsilon}(t,x)=\frac{1}{V(x)}\int_{\partial
D_{x}}\frac{1}{|\gamma_{2}^{\epsilon}(x,\epsilon y)|}c(x,\epsilon y,
u^{\epsilon}(t,x,\epsilon
y))dS_{x}.\label{NewNonLinearityBeforeLimit}
\end{equation}

Proceeding similarly now to Lemma 2.1 and recalling
that $v$ satisfies (\ref{NeymmanProblemForParabolicProblemAndu}), we see that there is a constant
$K^{\epsilon}=K(\|\nabla_{x}v^{\epsilon}\|,\|\nabla_{y}v^{\epsilon}\|,\|\triangle_{x}v^{\epsilon}\|,\|\triangle_{y}v^{\epsilon}\|,\|v^{\epsilon}_{t}\|,\|\gamma_{1}^{1}\|,T)$
such that:
\begin{eqnarray}
& &\sup_{0\leq t \leq T}E|\int_{0}^{t}\frac{1}{2}\overline{c}^{\epsilon}(t-s,X^{\epsilon}_{s})ds-
\int_{0}^{t}\epsilon
c(X^{\epsilon}_{s},Y^{\epsilon}_{s},u^{\epsilon}(t-s,X^{\epsilon}_{s},Y^{\epsilon}_{s}))dL^{\epsilon}_{s} |^{2}\nonumber\\
&\leq& \epsilon^{2}K^{\epsilon}(1+\sup_{0\leq t \leq T}E[\epsilon
L^{\epsilon}_{t}]^{2})\label{LimitIntegralWithLocalTimeAndFunction1}
\end{eqnarray}
We observe that $K^{\epsilon}$ depends on $\epsilon$ only through
functions that are uniformly bounded in $\epsilon$ (Proposition 3.7).
This observation and Lemma 2.2 imply that as $\epsilon\rightarrow
0$:
\begin{equation}
\sup_{0\leq t \leq T}E|\int_{0}^{t}\frac{1}{2}\overline{c}^{\epsilon}(t-s,X^{\epsilon}_{s})ds-
\int_{0}^{t}\epsilon
c(X^{\epsilon}_{s},Y^{\epsilon}_{s},u^{\epsilon}(t-s,X^{\epsilon}_{s},Y^{\epsilon}_{s}))dL^{\epsilon}_{s} |^{2}
\rightarrow 0. \label{LimitIntegralWithLocalTimeAndFunction2}
\end{equation}
Moreover the Lebesque dominated convergence Theorem, Lemma 3.3, the compactness of the family $\{u^{\epsilon}\}$ and (\ref{NewNonLinearityBeforeLimit}), imply that as $\epsilon\rightarrow 0$:
\begin{equation}
\sup_{0\leq t \leq
T}E|\int_{0}^{t}\frac{1}{2}\overline{c}^{\epsilon}(t-s,X^{\epsilon}_{s})ds-
\int_{0}^{t}\overline{c}(t-s,X_{s})ds |^{2} \rightarrow 0,
\label{LimitIntegralWithLocalTimeAndFunction3}
\end{equation}
where $\overline{c}^{\epsilon}, \overline{c}$ and $X_{t}$ are given by (\ref{NewNonLinearityBeforeLimit}), (\ref{NewNonLinearity}) and (\ref{LimitingStochasticProcess}) respectively.

Now let $u^{\epsilon}(t,x,y), u(t,x)$ be the solutions to (\ref{FeynmanKacFormulaBeforeLimit}) and (\ref{FeynmanKacFormulaAfterLimit}) respectively.
Taking into account relations (\ref{LimitIntegralWithLocalTimeAndFunction2}), (\ref{LimitIntegralWithLocalTimeAndFunction3}), the weak convergence of $X^{\epsilon}_{t}$ to $X_{t}$ as $\epsilon\rightarrow 0$ (which is implied by Lemma 3.3) and Proposition 3.2 we get the statement of the Theorem.
\begin{flushright}
$\square$
\end{flushright}

We conclude this section with the a-priori bounds for the H\"{o}lder
norm of the solution and for the sup-norm of the solution, the first
and the second derivatives of the solution of
(\ref{InitialEquationWith_e}). These bounds will be uniform in
$\epsilon$. The method follows closely \cite{AF1}.

\vspace{0.2cm}

Let us first introduce some notation.

We write $U_{T}^{\epsilon}=[0,T) \times D^{\epsilon}$, $\overline{U}_{T}^{\epsilon}=[0,T) \times \overline{D}^{\epsilon}$, $\partial U_{T}^{\epsilon}=[0,T) \times \partial D^{\epsilon}$ and $V_{T}^{\epsilon}=(0,T) \times D^{\epsilon}$, where $\overline{D}^{\epsilon}=D^{\epsilon}\cup \partial D^{\epsilon}$.

For $0<a<1$, $T>0$ and for any function $g$ we write:
\begin{eqnarray}
\parallel g \parallel_{U_{T}^{\epsilon}}&=&\sup_{(t,z) \in U_{T}^{\epsilon}}|g(t,z)|\nonumber\\
\|H^{a}g\|_{U_{T}^{\epsilon}}&=&\sup_{(t,z),(t',z') \in U_{T}^{\epsilon}}\frac{|g(t,z)-g(t',z')|}{|t-t'|^{a/2}+|z-z'|^{a}}\nonumber\\
\|g\|_{U_{T}^{\epsilon},a}&=&\|g\|_{U_{T}^{\epsilon}}+\|H^{a}g\|_{U_{T}^{\epsilon}}\nonumber\\
\|g\|_{D^{\epsilon},T,1+a}&=&\|g\|_{U_{T}^{\epsilon},a}+\|g_{t}\|_{U_{T}^{\epsilon}}+\|Dg\|_{(0,T)\times D^{\epsilon}}\nonumber\\
\overline{\|g\|}_{D^{\epsilon},T,1+a}&=&\|g\|_{\overline{D}^{\epsilon},T,1+a}\nonumber
\end{eqnarray}

\vspace{0.2cm}

Moreover for notational convenience we will write $z=(x,y)$.

\vspace{0.2cm}

\textbf{Lemma 3.5.} Under our assumptions there exists a constant $C_{1}$, independent of $\epsilon>0$, such that
\begin{center}
$0\leq u^{\epsilon}\leq C_{1}$ in $\overline{U}_{T}^{\epsilon}$.
\end{center}

\textbf{Proof.} Lemma 3.5 can be proven using equation (\ref{FeynmanKacFormulaBeforeLimit}). Here we give an analytic proof of the claim.
For any fixed $b>0$ we define the function
\begin{displaymath}
w^{\epsilon}=(u^{\epsilon}-b)^{+}=\max \lbrace u^{\epsilon}-b,0\rbrace.
\end{displaymath}
It is easy to show that
\begin{displaymath}
w^{\epsilon}_{t} \leq \frac{1}{2}\triangle w^{\epsilon} \hspace{1cm} \textrm{on} \hspace{0.2cm} (0,T)\times D^{\epsilon}.
\end{displaymath}
in the weak sense. Let us choose now $b=\max\{N,\|f\|\}$, where $N$
is such that if $u>N$ then $c(\cdot,\cdot,u)<-M$ for some $M>0$. Then
$$w^{\epsilon}(0,x,y)=0.$$

Let us now assume that $w^{\epsilon}$ attains a maximum positive
value on the boundary $\partial V_{T}^{\epsilon}$ at the point
$(t_{o},x_{o},y_{o})$. Since $u^{\epsilon}$ is continuous up to the
boundary, there exists a connected set $\Delta$ such that
$(t_{o},x_{o},y_{o})\in \Delta$, $\Delta \subset \partial
V_{T}^{\epsilon}$ and $w^{\epsilon}>0$ on $\Delta$, i.e.
$u^{\epsilon}>b$ on $\Delta$. Since $(t_{o},x_{o},y_{o})$ is a
maximum for $w^{\epsilon}$ and $\gamma^{\epsilon}$ is the inward
normal derivative we get that $\frac{\partial w^{\epsilon}}{\partial
\gamma^{\epsilon}}\leq 0$ at $(t_{o},x_{o},y_{o})$. But on $\Delta$
we have that $\frac{\partial w^{\epsilon}}{\partial
\gamma^{\epsilon}}=\frac{\partial u^{\epsilon}}{\partial
\gamma^{\epsilon}}=-\epsilon c(x,y,u^{\epsilon})u^{\epsilon}$.
Taking into account the particular choice of $b$ and that $c(\cdot,\cdot,u)<-M$ for $u>N$, we
get that $-\epsilon c(x,y,u^{\epsilon})u^{\epsilon}>0$ at
$(t_{o},x_{o},y_{o})$. Thus we have a contradiction and so maximum
principle implies that
\begin{displaymath}
w^{\epsilon}=0\Longrightarrow u^{\epsilon}< b \hspace{0.2cm}
\textrm{in} \hspace{0.2cm} \overline{U}_{T}^{\epsilon}.
\end{displaymath}
Lastly maximum principle again implies that $u^{\epsilon}\geq 0$.
\begin{flushright}
$\square$
\end{flushright}
Let us consider the following linear parabolic pde:
\begin{eqnarray}
v^{\epsilon}_{t}&=&\frac{1}{2}\triangle v^{\epsilon}, \hspace{3.1cm} \textrm{in}\hspace{0.1cm} (0,T) \times D^{\epsilon}\label{InitialLinearEquationWith_e}\\
v^{\epsilon}(0,x,y)&=&f(x), \hspace{3.3cm} \textrm{on}\hspace{0.1cm} \{0\} \times D^{\epsilon}\nonumber\\
\frac{\partial v^{\epsilon}}{\partial \gamma^{\epsilon}}&=&-\epsilon
c(x,y)v^{\epsilon}, \hspace{2.3cm} \textrm{on} \hspace{0.1cm}(0,T)
\times \partial D^{\epsilon},\nonumber
\end{eqnarray}
where $f,c$ are bounded smooth functions. Under the standard hypotheses problem (\ref{InitialLinearEquationWith_e})
has a unique classical solution (Theorem 5.3.2 in \cite{AF1}).

\vspace{0.2cm}

\textbf{Lemma 3.6.} There is a constant $C$, independent of $\epsilon$, and an open set $I\subset (0,1)$ such that for any $a \in I$:
\begin{equation}
\overline{\|v^{\epsilon}\|}_{D^{\epsilon},T,1+a}+\|D^{2}v^{\epsilon}\|_{V_{T}^{\epsilon}}\leq C. \label{BoundForInitialLinearEquationWith_e}
\end{equation}
\textbf{Proof.}  We will give just a sketch of the proof, since the
analysis follows \cite{P1}, \cite{P2}, \cite{P3} and \cite{AF1}. The
calculations are lengthy but standard.

We solve the second initial-boundary value problem (\ref{InitialLinearEquationWith_e}) by reducing it to an integral equation, i.e. we write:
\begin{equation}
v^{\epsilon}(t,z)=\int_{0}^{t}\int_{\partial D^{\epsilon}}\Gamma^{\epsilon}(t,z,\tau,\xi)\phi^{\epsilon}(\tau,\xi)d\partial D^{\epsilon}_{\xi}d\tau +\int_{D^{\epsilon}}\Gamma^{\epsilon}(t,z,0,\xi)f(\xi)d\xi \label{Represenataion1},
\end{equation}
where $\Gamma^{\epsilon}(t,z,\tau,\xi)=(2\sqrt{\pi})^{-n-m}(t-\tau)^{-\frac{n+m}{2}}\exp[-\frac{(z_{i}-\xi_{i})(z_{j}-\xi_{j})}{4(t-\tau)}]$
is the fundamental solution to the heat equation and $\phi(t,z)$ is the solution to a Voltera type integral equation:
\begin{eqnarray}
\phi^{\epsilon}(t,z)&=&2\int_{0}^{t}\int_{\partial D^{\epsilon}}[\frac{\partial \Gamma^{\epsilon}(t,z,\tau,\xi)}{\partial \gamma^{\epsilon}}+ \epsilon c(z)\Gamma^{\epsilon}(t,z,\tau,\xi)]\phi^{\epsilon}(\tau,\xi)d\partial D^{\epsilon}_{\xi}d\tau  \label{volteraForPhi}\\
&+& 2[\int_{D^{\epsilon}}\frac{\partial
\Gamma^{\epsilon}(t,z,0,\xi)}{\partial \gamma^{\epsilon}}f(\xi)d\xi+
\epsilon
c(z)\int_{D^{\epsilon}}\Gamma^{\epsilon}(t,z,0,\xi)f(\xi)d\xi]\nonumber
\end{eqnarray}

Let us now define
\begin{eqnarray}
F^{\epsilon}(t,z)&=&\int_{D^{\epsilon}}\frac{\partial \Gamma^{\epsilon}(t,z,0,\xi)}{\partial \gamma^{\epsilon}}f(\xi)d\xi+ \epsilon c(z)\int_{D^{\epsilon}}\Gamma^{\epsilon}(t,z,0,\xi)f(\xi)d\xi \nonumber\\
M_{1}(t,z,\tau,\xi)&=&\frac{\partial \Gamma^{\epsilon}(t,z,\tau,\xi)}{\partial \gamma^{\epsilon}}+ \epsilon c(z)\Gamma^{\epsilon}(t,z,\tau,\xi)\nonumber\\
M_{\nu+1}(t,z,\tau,\xi)&=&\int_{0}^{t}\int_{\partial D^{\epsilon}}M_{1}(t,z,t',z')M_{\nu}(t',z',\tau,\xi)d\partial D^{\epsilon}_{z'}dt'\nonumber
\end{eqnarray}

It can be shown (see \cite{AF1}) that there is a H\"{o}lder continuous (in space variables) and bounded (with bound and H\"{o}lder coefficient independent of $\epsilon$) solution $\phi^{\epsilon}$ for (\ref{volteraForPhi}), expressed in the form:
\begin{equation}
\phi^{\epsilon}(t,x)=2F^{\epsilon}(t,z)+2 \sum_{\nu=1}^{\infty}\int_{0}^{t}\int_{\partial D^{\epsilon}}M_{\nu}(t,z,\tau,\xi)F^{\epsilon}(\tau,\xi)d\partial D^{\epsilon}_{\xi}d\tau \label{FinalPhi}
\end{equation}
Using the boundedness and the H\"{o}lder continuity of
(\ref{FinalPhi}) and (\ref{Represenataion1}), one can show (see
\cite{P1}, \cite{P2}, \cite{P3} and \cite{AF1}) that there is a constant $C$,
independent of $\epsilon$, such that
$$\overline{\|v^{\epsilon}\|}_{D^{\epsilon},T,1+a}+\|D^{2}v^{\epsilon}\|_{V_{T}^{\epsilon}}\leq C.$$

\begin{flushright}
$\square$
\end{flushright}

Now we are ready to prove the result for the a-priori bounds:

\vspace{0.2cm}

\textbf{Proposition 3.7.} There is a constant $C$, independent of $\epsilon$, and an open set $I\subset (0,1)$ such that for any $b>a \in I$ ($a$ is the constant from Lemma 3.6.):
\begin{equation}
\overline{\|u^{\epsilon}\|}_{D^{\epsilon},T,1+b}+\|D^{2}u^{\epsilon}\|_{V_{T}^{\epsilon}}\leq C. \label{BoundForInitialEquationWith_e}
\end{equation}
where $u^{\epsilon}$ is a classical solution to (\ref{InitialEquationWith_e}).

\textbf{Proof.}
We will use Schauder's fixed point Theorem. Let us first define for convenience $\overline{\|\cdot\|}_{2+a}=\overline{\|\cdot\|}_{D^{\epsilon},T,1+a}+\|D^{2}\cdot\|_{V_{T}^{\epsilon}}$.

Let $\mathcal{C}^{2+a}$  be the Banach space of all functions
$u^{\epsilon}(t,z)$ that are continuous in
$\overline{U}_{T}^{\epsilon}$ with norm
$\overline{\|u^{\epsilon}\|}_{2+a}$.

For any $C>0$, let $\mathcal{C}_{C}^{2+a}$ be the set
$\{u^{\epsilon}:u^{\epsilon}\in \mathcal{C}^{2+a},
\overline{\|u^{\epsilon}\|}_{2+a}\leq C\}$.

For every $u^{\epsilon}\in \mathcal{C}_{C}^{2+a}$  define
$w^{\epsilon}=Tu^{\epsilon}$ to be the solution to the following
problem:

\begin{eqnarray}
w^{\epsilon}_{t}&=&\frac{1}{2}\triangle w^{\epsilon}, \hspace{2cm} \textrm{in}\hspace{0.1cm} (0,T) \times D^{\epsilon}\label{SchauderInitialEquationWith_e}\\
w^{\epsilon}(0,z)&=&f(x), \hspace{2.2cm} \textrm{on}\hspace{0.1cm} \{0\} \times D^{\epsilon}\nonumber\\
\frac{\partial w^{\epsilon}}{\partial \gamma^{\epsilon}}&=&-\epsilon
c(z,u^{\epsilon})w^{\epsilon}, \hspace{1cm} \textrm{on}
\hspace{0.1cm}(0,T) \times \partial D^{\epsilon},\nonumber
\end{eqnarray}

Then similarly as in Lemma 3.6, one can write:

\begin{equation}
w^{\epsilon}(t,z)=\int_{0}^{t}\int_{\partial D^{\epsilon}}\Gamma^{\epsilon}(t,z,\tau,\xi)\phi^{\epsilon}(\tau,\xi)d\partial D^{\epsilon}_{\xi}d\tau +\int_{D^{\epsilon}}\Gamma^{\epsilon}(t,z,0,\xi)f(\xi)d\xi,\label{SchauderInitialSolutionWith_e}
\end{equation}
where $\phi^{\epsilon}(t,z)$ satisfies:
\begin{eqnarray}
\phi^{\epsilon}(t,z)&=&2\int_{0}^{t}\int_{\partial D^{\epsilon}}[\frac{\partial \Gamma^{\epsilon}(t,z,\tau,\xi)}{\partial \gamma^{\epsilon}}+ \epsilon c(z,u^{\epsilon})\Gamma^{\epsilon}(t,z,\tau,\xi)]\phi^{\epsilon}(\tau,\xi)d\partial D^{\epsilon}_{\xi}d\tau  \nonumber\\
&+& 2[\int_{D^{\epsilon}}\frac{\partial
\Gamma^{\epsilon}(t,z,0,\xi)}{\partial \gamma^{\epsilon}}f(\xi)d\xi+
\epsilon
c(z,u^{\epsilon})\int_{D^{\epsilon}}\Gamma^{\epsilon}(t,z,0,\xi)f(\xi)d\xi]\label{SchauderVolteraForPhi}
\end{eqnarray}
We shall prove that $T$ has a fixed point.

Since $u^{\epsilon}$ and $c$ are bounded functions, one can show, in
the same way as in the proof of Lemma 3.6, that the function
$\phi^{\epsilon}(t,z)$ that satisfies (\ref{SchauderVolteraForPhi})
is bounded and H\"{o}lder continuous (in space variables) with bound
and H\"{o}lder constant independent of $\epsilon$.

Using this result and representation (\ref{SchauderInitialSolutionWith_e}) one can conclude (Lemma 3.6) that there is a constant $C$ such that
\begin{displaymath}
\overline{\|w^{\epsilon}\|}_{2+a}\leq C.
\end{displaymath}
So $T$ maps $\mathcal{C}_{C}^{2+a}$ into itself for an appropriately
chosen constant $C$.

Now let $\{u^{\epsilon}_{n}\}$ be a sequence of functions that
belong to $\mathcal{C}_{C}^{2+a}$ and
$w^{\epsilon}_{n},\phi^{\epsilon}_{n}$ be defined by
(\ref{SchauderInitialSolutionWith_e}) and
(\ref{SchauderVolteraForPhi}) when $u^{\epsilon}=u^{\epsilon}_{n}$.
Assume that
$\overline{\|u^{\epsilon}_{n}-u^{\epsilon}\|}_{2+a}\rightarrow
0$ as $n\rightarrow\infty$. We need to show that
$\overline{\|w^{\epsilon}_{n}-w^{\epsilon}\|}_{2+a}\rightarrow
0$ as $n\rightarrow\infty$.

The continuity of the function $c(z,u)$ in $u$-variables imply that
$\|\phi^{\epsilon}_{n}-\phi^{\epsilon}\|_{U_{T}^{\epsilon}}\rightarrow
0$ as $n\rightarrow\infty$. This and
(\ref{SchauderInitialSolutionWith_e}) give us
$\overline{\|w^{\epsilon}_{n}-w^{\epsilon}\|}_{2+a}\rightarrow
0$.

Therefore $T$ is a continuous map.

Next we need to show that $T$ maps $\mathcal{C}_{C}^{2+a}$ into a
compact subset of $\mathcal{C}_{C}^{2+a}$. This is an easy
consequence of Theorem 7.1.1 of \cite{AF1}, which states that for
$0<a<b<1$, the bounded subsets of  $\mathcal{C}^{2+b}$ are
pre-compact subsets of $\mathcal{C}^{2+a}$.

Lastly $\mathcal{C}_{C}^{2+b}$ is a closed convex set of the Banach
space $\mathcal{C}^{2+b}$.

Therefore by Schauder's Fixed Point Theorem we get that $T$ has a fixed point, i.e. there exists a $u^{\epsilon}$ such that $u^{\epsilon}=Tu^{\epsilon}$ and actually
$$u^{\epsilon}=Tu^{\epsilon} \in \mathcal{C}_{C}^{2+b}.$$

\begin{flushright}
$\square$
\end{flushright}

\section{Some Results On Wave Front Propagation}

In this section we will see some
applications of Theorem 3.4 to the question of wave front propagation
in narrow domains. As we mentioned in the introduction, corresponding results on the standard reaction diffusion equation (\ref{InitialEquationWithout_e}) (see chapter $6$ and $7$ in \cite{K1}, \cite{G1} and \cite{NX}) allow to describe the asymptotic wavefront motion for (\ref{InitialEquationWith_e}).

We will focus on two different cases. In subsection 4.1 we consider
the case where the functions $c(\cdot,u)$, $V(\cdot)$, $S(\cdot)$
and $f(\cdot)$ change slowly in $x$, i.e. $c(\cdot,u)=c(\delta
x,u)$, $V(\cdot)=V(\delta x)$, $S(\cdot)=S(\delta x)$ and
$f(\cdot)=f(\delta x)$ for $0<\delta\ll 1$. We first assume that the
nonlinear boundary term in (\ref{InitialEquationWith_e}),
$c(x,y,u)$, is of K-P-P type for $y=0$, i.e. $c(x,0,u)$ is positive
for $u<1$, negative for $u>1$ and $c(x)=c(x,0,0)=\max_{0\leq u\leq
1}c(x,0,u)$. We will see how the motion of the wavefront depends on
the behavior of the cross-sections $D_{x}$ of the domain $D$. In
particular, using the results of \cite{K1} (chapter 6) we will see
that in the case of the nonlinear term of K-P-P type and for $x\in
\mathbb{R}$ the wavefront can have jumps. Actually, the jumps of the
wavefront appear at positions where the tube becomes thinner. The
results are given in Theorem 4.1.2, Theorem 4.1.5 and Theorem
4.1.7. Then we briefly discuss the bistable case, i.e. when
$c(x,0,u)>0$ for $u\in(\mu,1)$ and $c(x,0,u)<0$ for $u\in
(0,\mu)\cup (1,\infty)$, where $0<\mu<1$. In this case we consider a
specific example and we will see how the asymptotic speed of the
wavefront depends on the surface area to volume ratio
$\frac{S(x)}{V(x)}$. In subsection 4.2, we return to the K-P-P case,
but now we consider front propagation when $x\in \mathbb{R}$ and the
boundary $\partial D^{1}$ of $D^{1}$ is determined by stationary
random processes on $\mathbb{R}$ on some probability space
$(\hat{\Omega},\hat{\mathfrak{F}},\hat{P})$. The conclusion is in
Theorem 4.2.7.

We will denote by
$\overline{c}(x,u):=\frac{1}{2}\frac{S(x)}{V(x)}c(x,0,u(t,x))$ the
nonlinear term in (\ref{InitialEquationWithout_e}). Obviously the
type of $\overline{c}(x,u)$  (K-P-P or bistable) is determined by
$c(x,0,u)$ and vice-versa.

\subsection{Wave Fronts in Slowly Changing Media}
\vspace{0.2cm}

Let us assume that the functions $c(\cdot,u)$, $V(\cdot)$, $S(\cdot)$ and $f(\cdot)$ change slowly in $x$, i.e.
$c(\cdot,u)=c(\delta x,u)$, $V(\cdot)=V(\delta x)$, $S(\cdot)=S(\delta x)$ and $f(\cdot)=f(\delta x)$ for $0<\delta\ll 1$.

We start with the case where the nonlinear term $\overline{c}(x,u)$
of (\ref{InitialEquationWithout_e}) is of K-P-P type. We
additionally assume that the closure of the support of $f$, $F_{o}$,
coincides with the closure of its interior. Lastly we take for
brevity $x\in \mathbb{R}^{1}$ and
$\overline{c}(x)=\overline{c}(x,0)=\frac{1}{2}\frac{S(x)}{V(x)}c(x,0,0)$
(recall that $c(x,0,0)=\sup_{0\leq u\leq 1}c(x,0,u)$) to be an
increasing function.

\vspace{0.2cm}

Let $\phi :[0,T]\rightarrow
\mathbb{R}^{1}$ and introduce the functional
\begin{equation}
R_{0,T}(\phi)=\cases{\int_{0}^{T}[\overline{c}(\phi_{s})-\frac{1}{2}|\dot{\phi}_{s}|^{2}]ds,\label{R_functional} & $\phi$ is absolutely continuous \cr
                                                     +\infty, & for the rest of  $\mathcal{C}_{0,T}.$\cr  }
\end{equation}
Put
\begin{equation}
W(t,x)=\sup\{R_{0,t}(\phi):\phi \in \mathcal{C}_{0,t}(\mathbb{R}^{1}), \phi_{0}=x, \phi_{t}\in F_{o}\}.\label{DefinitionOfW}
\end{equation}
We say that condition (N) is satisfied if for any $t>0$ and $(t,x)\in
\{(t,x):W(t,x)=0\}:$

$$W(t,x)=\sup\{
R_{0,t}(\phi):\phi_{0}=x, \phi_{t}\in
F_{o}, (t-s,\phi_{s})\in \{(t,x):W(t,x)<0\}\}.$$

As it is mentioned in chapter 10 of \cite{FW1}, condition (N) is fulfilled for the smooth and increasing function $\overline{c}(x)$. Moreover as we shall see in Theorem 4.1.2, $W(t,x)$ determines the motion of the wave
front for $u^{\epsilon}$ for small enough $\epsilon>0$.

Let us consider $u(t,x)$, the solution to equation (\ref{InitialEquationWithout_e}), for $n=1$. If we set
$u^{\delta}(t,x)=u(t/\delta,x/\delta)$, then $u^{\delta}$ is the
solution to the following parabolic problem:
\begin{eqnarray}
u^{\delta}_{t}&=&\frac{\delta}{2} u_{xx}^{\delta} +\frac{\delta}{2} \frac{ V_{x}(x)}{ V(x)} u_{x}^{\delta} + \frac{1}{\delta }\overline{c}(x,u^{\delta}(t,x))u^{\delta}, \hspace{0.1cm} \textrm{in}\hspace{0.1cm} (0,\infty) \times \mathbb{R}^{1}\nonumber\\
u^{\delta}(0,x)&=&f(x)\geq 0, \hspace{0.2cm} \textrm{on}\hspace{0.1cm}
\{0\} \times \mathbb{R}^{1}.\label{FinalEquationWith_d}
\end{eqnarray}

Under the assumptions above, the following theorem, which is a reformulation of Theorem 6.2.1 of \cite{K1}, states that $W(t,x)$ determines the motion of the wave front for $u^{\delta}(t,x)$ under condition (N):

\vspace{0.2cm}

\textbf{Theorem 4.1.1.} Let $u^{\delta}(t,x)$ be the solution to (\ref{FinalEquationWith_d}). Then under condition (N) we have:
\begin{eqnarray}
\lim_{\delta \downarrow 0}u^{\delta}(t,x)=\cases{1, & $W(t,x)>0$ \cr
                                                     0, & $W(t,x)<0.$ \cr  }
\end{eqnarray}

Let us consider now equation (\ref{InitialEquationWith_e}) for $n=1$, $c(\cdot,u)=c(\delta x,u)$, $f(\cdot)=f(\delta x)$ in a slowly changing in $x$ narrow domain $D^{\epsilon,\delta}$, so that $V(\cdot)=V(\delta x)$, $S(\cdot)=S(\delta x)$. Let us define $u^{\epsilon,\delta}(t,x,y)=u^{\epsilon}(t/\delta,x/\delta,y)$. Under the assumptions above, Theorems 3.4 and 4.1.1 imply that $W(t,x)$ will determine the motion of the wave
front in this case too, as follows:

\vspace{0.2cm}

\textbf{Theorem 4.1.2.} The following statement holds:
\begin{eqnarray}
\lim_{\delta\downarrow 0}\lim_{\epsilon\downarrow 0} u^{\epsilon,\delta}(t,x,y)=\cases{1, & $W(t,x)>0$ \cr
                                                     0, & $W(t,x)<0.$ \cr  }
\end{eqnarray}

So the equation $W(t,x)=0$ defines the position of the interface (wavefront) between areas where $u^{\epsilon,\delta}$ (for $\epsilon>0$ and $\delta>0$ small enough) is close to $0$ and to $1$. Actually, as we shall see below the wavefront may have jumps.
It is known (see
chapter 6 in \cite{K1}),
that because of the dependance of $\overline{c}(x)$ on $x$, the wave
front of $u^{\delta}$ may have jumps and new sources may be
"igniting" ahead of the front. We will give sufficient conditions that guarantee such jumps for a class of smooth and increasing functions $\bar{c}(x)$. Hence Theorem 4.1.2 implies that one can predict appearances of
new sources and jumps of the wave front of $u^{\epsilon,\delta}$ for
$\epsilon>0$ and $\delta>0$ small enough.
Reaction-Diffusion Eqautions (RDE's) with Nonlinear Boundary Conditions In Narrow Domains
Let $t^{*}=t^{*}(x,\bar{c}(\cdot))$ be such that $W(t^{*},x)=0$. Such a $t^{*}(x,\bar{c}(\cdot))$ is defined in a unique way.
\newpage

\begin{figure}[ht]
\begin{center}
\includegraphics[scale=1, width=8 cm, height=4.8 cm]{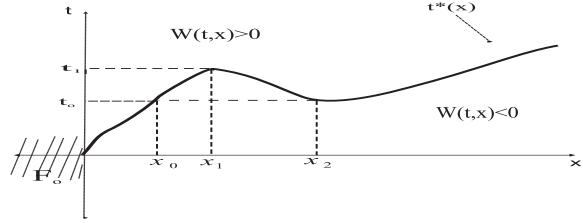}
\caption[A wavefront that has jumps]{ A wavefront that jumps from $x_{0}$ to $x_{2}$ at time $t_{0}$.}
\end{center}
\end{figure}

We have the following proposition (see chapter 6 in \cite{K1} for more details):

\vspace{0.2cm}

\textbf{Proposition 4.1.3.} Let $t^{*}(x)$ be as in Figure 1 and $F_{o}=\{x\in \mathbb{R}^{1},x<0\}$. Then the wavefront jumps from $x_{o}$ to $x_{2}$ at time $t_{o}$ (see Figure 1), i.e.:
\begin{enumerate}
\item{If $t\leq t_{0}$ then  $\lim_{\delta\downarrow 0}\lim_{\epsilon\downarrow 0} u^{\epsilon,\delta}(t,x,y)=1$  for a connected set: \\
 $F_{t}=\lbrace x\in \mathbb{R}^{1}: W(t,x)>0 \textrm{ and } x<x_{0}\rbrace$.}
\item{If $t_{0}<t<t_{1}$ then the set where $\lim_{\delta\downarrow 0}\lim_{\epsilon\downarrow 0} u^{\epsilon,\delta}(t,x,y)=1$
consists of two connected components: \\
$F_{t}=\lbrace x\in \mathbb{R}^{1}: W(t,x)>0 \textrm{ and } x<x_{1}\rbrace\cup\lbrace x\in \mathbb{R}^{1}: W(t,x)>0 \textrm{ and }x>x_{1}\rbrace$.\\
The set $\lbrace x\in \mathbb{R}^{1}: W(t,x)>0 \textrm{ and } x<x_{1}\rbrace$ is at a positive distance from the set $\lbrace x\in \mathbb{R}^{1}: W(t,x)>0 \textrm{ and }x>x_{1}\rbrace$ for $t_{0}<t<t_{1}$.}
\item{If $t\geq t_{1}$ then  $\lim_{\delta\downarrow 0}\lim_{\epsilon\downarrow 0} u^{\epsilon,\delta}(t,x,y)=1$  for a connected set: \\  $F_{t}=\lbrace x\in \mathbb{R}^{1}: W(t,x)>0\rbrace$.}
\end{enumerate}

\vspace{0.2cm}

Based now on comparison results (Lemma 4.1.4) we will give sufficient conditions that guarantee jumps of the wavefront. In particular we will prove (Theorem 4.1.5) that if $\overline{c}(x)$ is a rapidly increasing smooth function, then  $t^{*}=t^{*}(x,\bar{c}(\cdot))$ such that $W(t^{*},x)=0$ is as in Figure  1.

The functional $R_{0,T}(\phi)$ defined in (\ref{R_functional}) and the function $W(t,x)$ defined in (\ref{DefinitionOfW}) depend also on $\bar{c}$. Hence we will write sometimes $R_{0,T}(\phi,\bar{c}(\cdot))$ and $W(t,x,\bar{c}(\cdot))$ in order to emphasize this dependence.

\vspace{0.1cm}

We have the following comparison result:

\vspace{0.2cm}

\textbf{Lemma 4.1.4.}
\begin{enumerate}
\item{Let $A$ be a positive number. Then $t^{*}(x,Ac(\cdot))=\frac{1}{\sqrt{A}}t^{*}(x,c(\cdot))$.}
\item{Let $a$ be a positive number and define $c_{a}(x)=c(ax)$. Then $t^{*}(x,c_{a}(\cdot))=\frac{1}{a}t^{*}(ax,c(\cdot))$.}
\item{Let $c_{1},c_{2}$ be two functions such that $c_{1}(x)<c_{2}(x)$ for every $x\in \mathbb{R}^{1}$. Then $t^{*}(x,c_{1}(\cdot))>t^{*}(x,c_{2}(\cdot))$.}
\end{enumerate}

\textbf{Proof.} Let us write $t^{*}_{A}=t^{*}(x,Ac(\cdot))$ and let $\phi^{A}$ be the extremal so that $W(t^{*}_{A},x,Ac(\cdot))=R_{0,t^{*}_{A}}(\phi^{A},Ac(\cdot))=0$. Such an extremal satisfies the following Euler-Lagrance equation:
\begin{eqnarray}
\ddot{\phi}^{A}(s)&=&-Ac'(\phi^{A}(s))\nonumber\\
\phi^{A}(0)&=& x\label{EquationForTheExtremal1}\\
\phi^{A}(t^{*}_{A})&=&0. \nonumber
\end{eqnarray}

Let us define now the function $\phi(s)=\phi^{A}(s/\sqrt{A})$. We claim that the function $\phi(s)$ is the extremal so that $W(\sqrt{A}t^{*}_{A},x,c(\cdot))=R_{0,\sqrt{A}t^{*}_{A}}(\phi,c(\cdot))=0$. Indeed it is easy to see that the definition of $\phi$ and the fact that  $R_{0,t^{*}_{A}}(\phi^{A},Ac(\cdot))=0$ imply that $R_{0,\sqrt{A}t^{*}_{A}}(\phi,c(\cdot))=0$. Moreover $\phi$ satisfy an Euler-Lagrange equation of the form (\ref{EquationForTheExtremal1}) with $Ac(x)$ and $t^{*}_{A}$ replaced by $c(x)$ and $\sqrt{A}t^{*}_{A}$ respectively. This proves the claim, which implies part (i) of the lemma.

\vspace{0.1cm}

Part (ii) of the lemma can be proven in a similar way. We define $t^{*}_{a}=t^{*}(x,c_{a}(\cdot))$ and let $\phi^{a}$ to be the extremal so that $W(t^{*}_{a},x,c_{a}(\cdot))=R_{0,t^{*}_{a}}(\phi^{a},c_{a}(\cdot))=0$. Then similarly as it is done in part (i), one should consider the function $\phi(s)$ that is defined by $\phi(s)=a\phi^{a}(s/a)$.

\vspace{0.05cm}

We prove now part (iii) of the lemma. Let us define $t^{*}_{1}=t^{*}(x,c_{1}(\cdot))$ and $t^{*}_{2}=t^{*}(x,c_{2}(\cdot))$. Moreover let $\phi^{1}$ be the extremal so that $W(t^{*}_{1},x,c_{1}(\cdot))=R_{0,t^{*}_{1}}(\phi^{1},c_{1}(\cdot))=0$. Since $c_{1}(x)<c_{2}(x)$ we have
\begin{equation}
0=R_{0,t^{*}_{1}}(\phi^{1},c_{1}(\cdot))<R_{0,t^{*}_{1}}(\phi^{1},c_{2}(\cdot)).\label{ComparisonLemma1}
\end{equation}
Furthermore, it is easy to see that $W(t,x)$ is an increasing function of $t$.

Let us assume now that $t^{*}_{1}\leq t^{*}_{2}$. This assumption and the fact that $W(t^{*}_{2},x,c_{2}(\cdot))=0$ imply that $W(t^{*}_{1},x,c_{2}(\cdot))\leq 0$. By recalling the definition of function $W$, one easily concludes that:
\begin{equation}
R_{0,t^{*}_{1}}(\phi^{1},c_{2}(\cdot))\leq 0.\label{ComparisonLemma2}
\end{equation}
However inequality (\ref{ComparisonLemma2}) contradicts (\ref{ComparisonLemma1}). Therefore $t^{*}(x,c_{1}(\cdot))>t^{*}(x,c_{2}(\cdot))$.
\begin{flushright}
$\square$
\end{flushright}
In section 6.2 of \cite{K1}, it is proven that if $\bar{c}(x)$,
instead of the smooth function
$\frac{1}{2}\frac{S(x)}{V(x)}c(x,0,0)$, is a piecewise constant
function, denoted by $d(x)$, such that
\begin{eqnarray}
d(x)=\cases{d_{1}, & $ x<x_{2} $ \cr
                                                     d_{2}, & $x\geq x_{2}.$ \label{StepFunction}  }
\end{eqnarray}
with $d_{2}>2 d_{1}>0$, then the function $t^{*}=t^{*}(x,d(\cdot))$ such that $W(t^{*},x,d(\cdot))=0$ is not monotone, as in Figure 1. More specifically the curves connecting the point $(0,0)$ with $(x_{1},t_{1})$ and $(x_{1},t_{1})$ with $(x_{2},t_{0})$ are line segments and for $x>x_{2}$,  $t^{*}=t^{*}(x,d(\cdot))$ is the solution to
$$\sup_{t}\lbrace d_{2}(t^{*}-t)+d_{1}t-\frac{(x-x_{2})^{2}}{2(t^{*}-t)}-\frac{x_{2}^{2}}{2t}\rbrace=0.$$
Moreover in this case
\begin{eqnarray}
t_{0}&=&x_{2}\frac{\sqrt{2(d_{2}-d_{1})}}{d_{2}} \label{TimesOfJump0}\\
t_{1}&=&\frac{1}{2\sqrt{2d_{1}}}(x_{2}+\sqrt{2d_{1}}t_{0}) \label{TimesOfJump1}
\end{eqnarray}
We will write $t_{0}=t_{0}(d)$ and $t_{1}=t_{1}(d)$ to emphasize the dependence of $t_{0}$ and $t_{1}$ on the function $d(x)$.

\vspace{0.05cm}

With the help of the result above and Lemma 4.1.4  we will give sufficient conditions
that guarantee jumps of the wavefront of $u^{\delta}(t,x)$ (and by Theorem 4.1.2 of $u^{\epsilon}(t,x,y)$ for $\epsilon>0$ and $\delta>0$ small enough) for a class of smooth and increasing functions $\bar{c}$.

\vspace{0.05cm}

Let us define the set
\begin{equation}
\Delta=\{(d_{1},d_{2})\in \mathbb{R}_{+}^{1}\times\mathbb{R}_{+}^{1}: d_{2}>2d_{1}\textrm{ and }d_{2}>2\sqrt{d_{1}(d_{2}-d_{1})}\}. \label{AdmissibleSet}
\end{equation}
It is easy to see that $\Delta$ is a non-empty set.

\vspace{0.1cm}

\textbf{Theorem 4.1.5.} Let $d(x)$ be the step function defined in (\ref{StepFunction}) such that $(d_{1},d_{2})\in \Delta$.
Consider real numbers $A$ and $a$ such that
\begin{enumerate}
\item{$a,A>1$.}
\item{$a\sqrt{A}<\frac{1}{2}[1+\frac{d_{2}}{2\sqrt{d_{1}(d_{2}-d_{1})}}]$.}
\end{enumerate}
Then for any smoothly increasing function $\bar{c}(x)$ such that
\begin{equation}
d(x)<\bar{c}(x)<Ad(ax)\label{SufficientConditionForAJump1}
\end{equation}
the wavefront corresponding to $\bar{c}$ has jumps. In particular the excitation reaches the region $\lbrace x>\frac{x_{1}}{a}+\delta \rbrace$ before it reaches the point $\frac{x_{1}}{a}$, where $\delta$ is a small enough positive number and $\frac{x_{1}}{a}$ is as in Figures 2 and 3.

\vspace{0.05cm}

\textbf{Proof.} Let us define $\bar{d}(x)=Ad(ax)$. Since $a,A>1$, the function d(x)
 is shifted vertically upwards and horizontally to the left.  So we get that $d(x)<\bar{d}(x)$ (see Figure
 2).

Parts (i) and (ii) of Lemma 4.1.4 imply that
$t^{*}(x,\bar{d}(\cdot))=\frac{1}{a\sqrt{A}}t^{*}(ax,d(\cdot))$. This and part
(iii) of Lemma 4.1.4 give that if $\bar{c}$ satisfies
(\ref{SufficientConditionForAJump1}), then $t^{*}(x,\bar{c}(\cdot))$ will
satisfy (see Figure 3):
\begin{equation}
\frac{1}{a\sqrt{A}}t^{*}(ax,d(\cdot))<t^{*}(x,\bar{c}(\cdot))<t^{*}(x,d(\cdot)).\label{RelationForWavefrontPosition1}
\end{equation}

\begin{figure}[hbp]
\begin{center}
\begin{minipage}[b]{.5\textwidth}
\centering
\includegraphics[scale=.4,width=6 cm, height=5 cm]{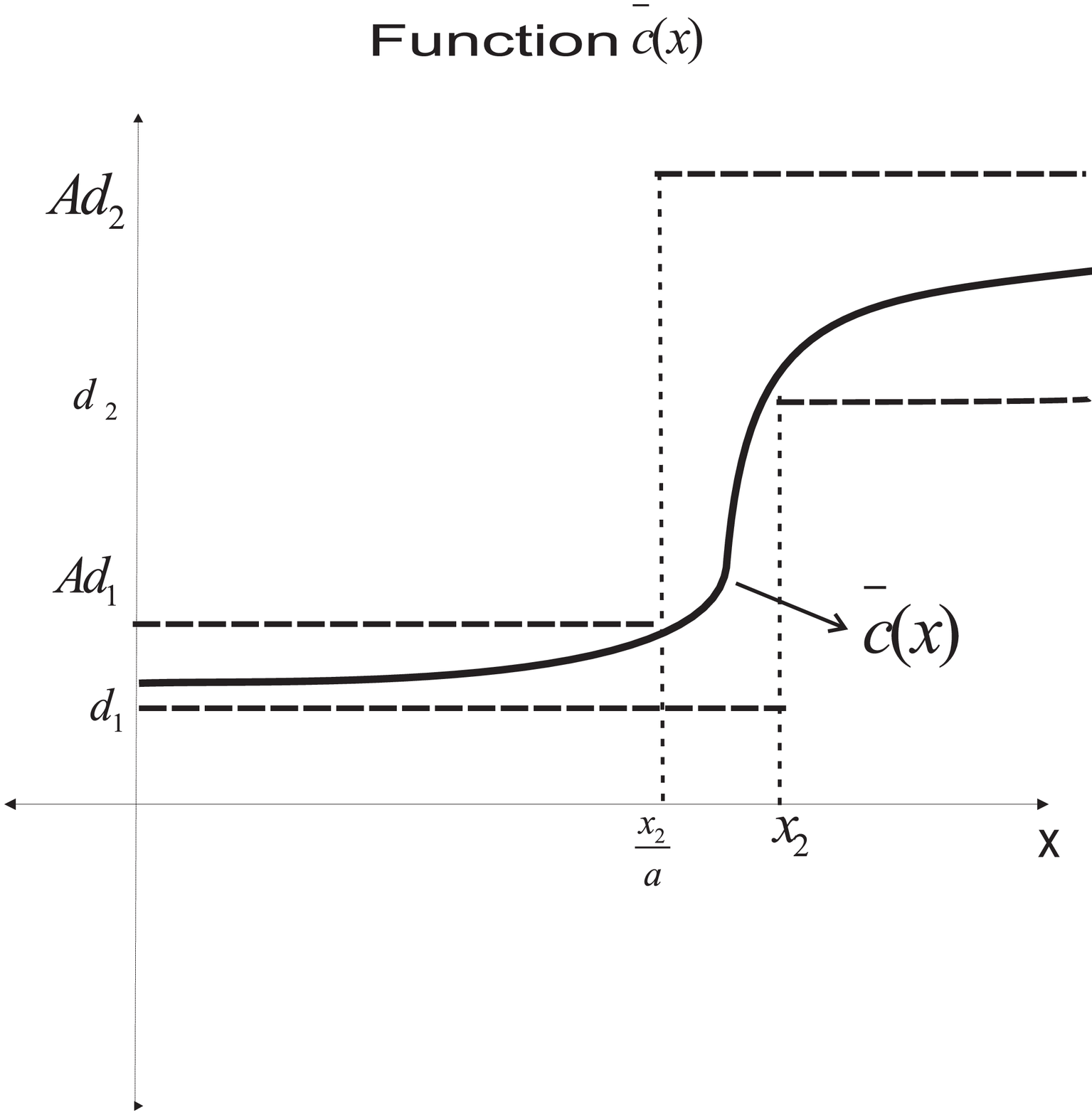}
\caption{$d(x)<\bar{c}(x)<Ad(ax)$}
\end{minipage}%
\begin{minipage}[b]{.5\textwidth}
\centering
\includegraphics[scale=.4,width=6 cm, height=5 cm]{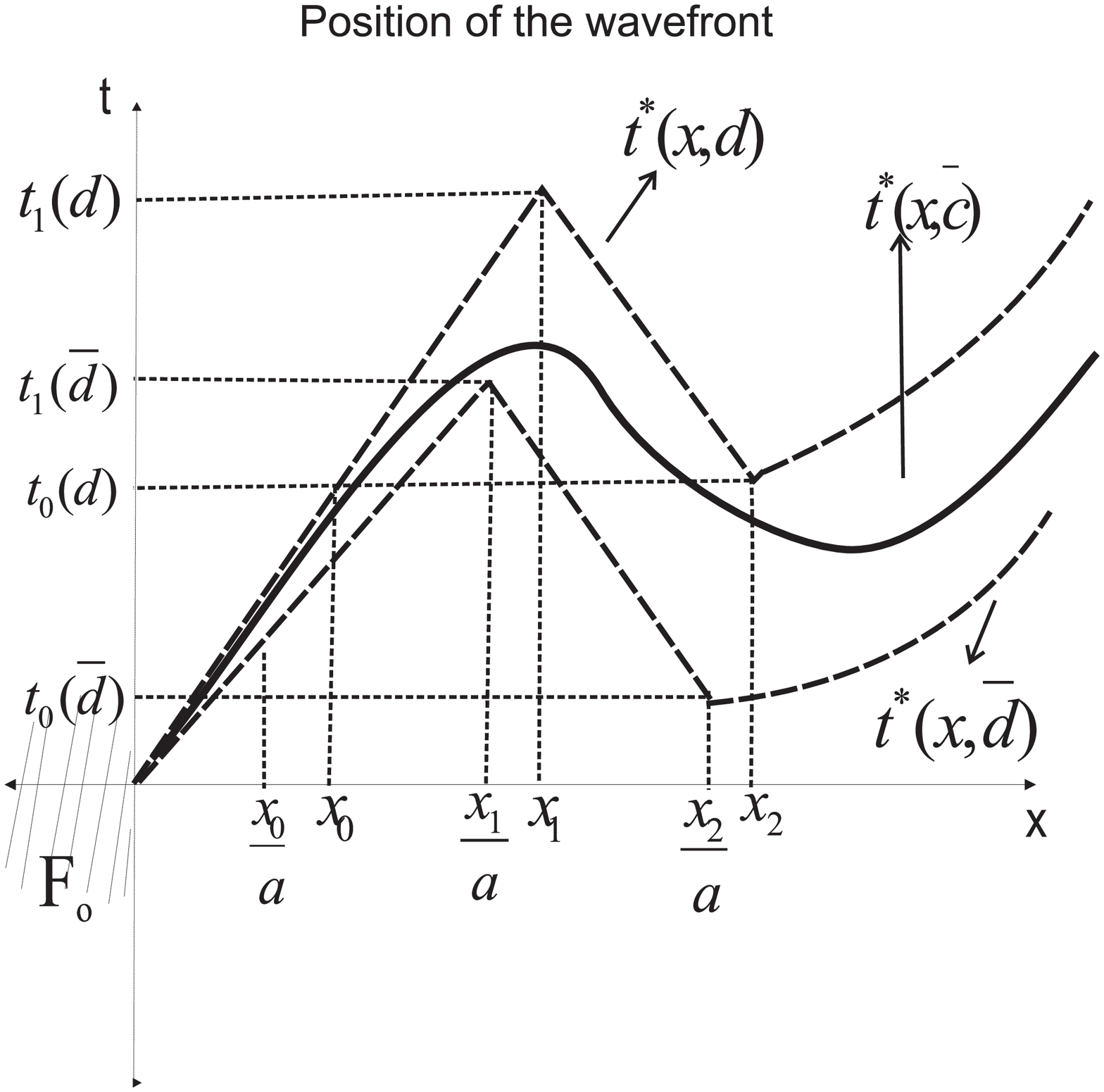}
\caption{$t^{*}(x,\bar{d})<t^{*}(x,\bar{c})<t^{*}(x,d)$}
\end{minipage}%
\end{center}
\end{figure}

We know that $t^{*}(x,d(\cdot))$ and $t^{*}(x,\bar{d}(\cdot))$ are not monotone (recall that $d$ and $\bar{d}$ are piecewise constant functions). We will show
that $t^{*}(x,\bar{c}(\cdot))$ is also not monotone (i.e it is as in Figure 1). Let us assume
that
\begin{equation}
t_{1}(\bar{d})>t_{0}(d), \label{SufficientConditionForAJump2}
\end{equation}
where $t_{0}(d)$ is as in (\ref{TimesOfJump0}) and $t_{1}(\bar{d})$
is defined similarly to $t_{1}(d)$ in (\ref{TimesOfJump1}) with
$d_{1},d_{2},x_{2}$ replaced by $Ad_{1},Ad_{2},\frac{x_{2}}{a}$
respectively. In particular (\ref{SufficientConditionForAJump2})
holds if condition (ii) above holds, i.e. if
$a\sqrt{A}<\frac{1}{2}[1+\frac{d_{2}}{2\sqrt{d_{1}(d_{2}-d_{1})}}]$.
Moreover the condition $d_{2}>2\sqrt{d_{1}(d_{2}-d_{1})}$ in the
definition of $\Delta$ in (\ref{AdmissibleSet}) implies that
$\frac{1}{2}[1+\frac{d_{2}}{2\sqrt{d_{1}(d_{2}-d_{1})}}]>1$, which
has to be true since $a,A>1$.

Inequality (\ref{SufficientConditionForAJump2}) can be equivalently written as $t^{*}(\frac{x_{1}}{a},\bar{d}(\cdot))>t(x_{2},d(\cdot))$. By this
and (\ref{RelationForWavefrontPosition1})  we immediately get that
\begin{equation}
t^{*}(x_{2},\bar{c}(\cdot)) < t^{*}(\frac{x_{1}}{a},\bar{c}(\cdot))
\end{equation}
which, since $\frac{x_{1}}{a}<x_{1}<x_{2}$, implies that
$t^{*}(x,\bar{c}(\cdot))$ is as in Figure 1 and so new sources are igniting
ahead of the wavefront.

In Figures 2 and 3 we see an illustration of the construction.
\begin{flushright}
$\square$
\end{flushright}
\textbf{Example.} An example of a function $\bar{c}(x)$ that satisfies the requirements of Theorem 4.1.5 is
\begin{equation}
\bar{c}(x)=\frac{Ad_{2}\mu+d_{1}e^{-\lambda(x-k)}}{\mu+e^{-\lambda(x-k)}},
\end{equation}
where $(d_{1},d_{2})\in \Delta$, $a,A$ satisfy assumptions (i) and (ii) of Theorem 4.1.5,  $k\in(\frac{x_{2}}{a},x_{2})$ and the constants $\mu$ and $\lambda$ are chosen so that $\bar{c}(\frac{x_{2}}{a})<Ad_{1}$ and $\bar{c}(x_{2})>d_{2}$.

In particular now if $\overline{c}(x)=\frac{1}{2}\frac{S(x)}{V(x)}$,
i.e. $c(x,0,0)=1$, is an increasing smooth function that satisfies
the requirements of Theorem 4.1.5, then the jump of the wavefront of
$u^{\epsilon,\delta}(t,x,y)$, for $\epsilon>0$ and $\delta>0$ small
enough, occurs when $\frac{S(x)}{V(x)}$ increases rapidly. This
implies, at least when the tube $D^{1}$ retains its shape as $x$
increases, that the jumps of the wave front occur at places where
the tube $D^{1}$ becomes thinner, i.e. when $V(x)$ decreases
significantly.

\vspace{0.05cm}

\textbf{Remark 4.1.6.}  Similar results hold for layers as well, i.e. for $x\in \mathbb{R}^{n}$ with $n>1$.

\vspace{0.05cm}

Using the results in \cite{F3} one can consider the limiting
behavior as $\delta,\epsilon \downarrow 0$ of
$u^{\epsilon,\delta}(t,x,y)$ when condition (N) is not fulfilled. We
will briefly discuss the result for the general case $x\in \mathbb{R}^{n}$.

Instead now of function $W(t,x)$ defined by (\ref{DefinitionOfW}), we consider the function
\begin{eqnarray}
W^{*}(t,x)=\sup \lbrace  \min_{0\leq s\leq t}R_{0,s}(\phi)&:&
\phi \in \mathcal{C}_{0,t}(\mathbb{R}^{n}) \textrm{ is absolutely continuous, }\nonumber \\
& & \phi_{0}=x, \hspace{0.2cm} \phi_{t}\in F_{o}\rbrace.\label{DefinitionOfW^{*}}
\end{eqnarray}

One can prove that $W^{*}(t,x)$ is Lipschitz continuous and that $W^{*}(t,x)\leq \min\{0, W(t,x)\}$.

Then Theorem 2.1 in \cite{F3} and Theorem 3.4 imply that $W^{*}(t,x)$ determines the motion of the wave front as follows:

\vspace{0.2cm}

\textbf{Theorem 4.1.7.} The following statements hold:
\begin{enumerate}
\item{For any compact subset $\Theta_{1}$ of the interior of $\lbrace (t,x): t>0,
W^{*}(t,x)=0\rbrace$,
\begin{displaymath} \lim_{\delta\downarrow
0}\lim_{\epsilon\downarrow 0} u^{\epsilon,\delta}(t,x,y)= 1 \textrm{
uniformly in } (t,x)\in \Theta_{1}.
\end{displaymath}}
\item{For any compact subset $\Theta_{2}$ of $\lbrace
(t,x):W^{*}(t,x)<0\rbrace$,
\begin{displaymath} \lim_{\delta\downarrow
0}\lim_{\epsilon\downarrow 0} u^{\epsilon,\delta}(t,x,y)= 0 \textrm{
uniformly in } (t,x)\in \Theta_{2}.
\end{displaymath}}
\end{enumerate}

\vspace{0.2cm}

We conclude subsection 4.1 with the case that the nonlinear term $\overline{c}(x,u)$
of (\ref{InitialEquationWithout_e}) is of bistable type, i.e. $\overline{c}(x,u)>0$ for $u\in(\mu,1)$, $\overline{c}(x,u)<0$ for
$u\in (0,\mu)\cup (1,\infty)$, where $0<\mu<1$. This problem was considered in
\cite{G1} and it was also presented in section 6.4 of \cite{K1}.

Here we restrict the analysis to a concrete example that allows to
give an exact formula for the asymptotic speed of the wavefront of
$u^{\epsilon,\delta}$ for $\epsilon>0$ and $\delta>0$ small enough.
As we will see the asymptotic speed of the wavefront is proportional to the square root of the
surface area to volume ratio $ \sqrt{\frac{S(x)}{V(x)}}$.

To be specific let
$x\in \mathbb{R}^{n}$, $c(x,0,u)=(u-\mu)(1-u)$, $0<\mu<\frac{1}{2}$
and assume that the function $u^{\delta}(t,x)$ (compare with
(\ref{FinalEquationWith_d})) is the solution to
\begin{eqnarray}
u^{\delta}_{t}&=&\frac{\delta}{2V(x)} \textrm{div}(V(x)\nabla_{x} u^{\delta}) + \frac{1}{\delta}\frac{1}{2}\frac{S(x)}{V(x)}(u^{\delta}-\mu)(1-u^{\delta})u^{\delta}, \hspace{0.3cm} \textrm{in}\hspace{0.1cm} (0,\infty) \times \mathbb{R}^{n}\nonumber \\
u^{\delta}(0,x)&=&f(x), \hspace{0.2cm} \textrm{on}\hspace{0.1cm} \{0\} \times \mathbb{R}^{n}.\label{FinalEquationWith_d1}
\end{eqnarray}
Consider a point $x\in \mathbb{R}^{n}$ to be excited at time $t$, if $u^{\delta}(t,x)$ (the solution to (\ref{FinalEquationWith_d1})) is close to $1$ and non-excited if  $u^{\delta}(t,x)$ is close to $0$. Then the Corollary of Theorem 4.1 of \cite{G1} gives us that for small $\delta>0$ the region $\lbrace x\in \mathbb{R}^{n}:f(x)>\mu \rbrace$ becomes excited and the region $\lbrace x\in \mathbb{R}^{n}:f(x)<\mu \rbrace$ becomes non-excited after a short starting phase. Now let $u^{\epsilon,\delta}(t,x,y)=u^{\epsilon}(t/\delta,x/\delta,y)$, where $u^{\epsilon}(t,x,y)$ is the solution to (\ref{InitialEquationWith_e}). Theorem 3.4 implies that the same conclusions hold for $u^{\epsilon,\delta}(t,x,y)$ for $\epsilon>0$ and $\delta>0$ small enough.

To compute the asymptotic propagation speed of excitation at $x\in \mathbb{R}^{n}$, let us consider the equation for the wave profile:
\begin{eqnarray}
& &\frac{1}{2}v_{\xi\xi}''(\xi)+a(x)v_{\xi}'(\xi)+ \frac{1}{2}\frac{S(x)}{V(x)}(v(\xi)-\mu)(1-v(\xi))v(\xi)=0, \hspace{0.2cm} \xi\in \mathbb{R}\label{RunningWaveBistableHomogeneousCaseExample}\\
& &\lim_{\xi\rightarrow -\infty}v(\xi)=1, \hspace{0.2cm}
\lim_{\xi\rightarrow \infty}v(\xi)=0.\nonumber
\end{eqnarray}
As it can be verified by direct substitution, equation (\ref{RunningWaveBistableHomogeneousCaseExample}) is
solvable if $a(x)$ is given by the formula
\begin{equation}
a(x)=\sqrt{\frac{1}{2}\frac{S(x)}{V(x)}}(\frac{1}{2}-\mu).\label{WaveFrontSpeed}
\end{equation}
Moreover, in our case, (\ref{WaveFrontSpeed}) is also the asymptotic propagation speed of excitation at $x\in \mathbb{R}^n$ and it is independent of direction.

Lastly, it is known that as the size of $D_{x}$ increases (without
changing shape), the surface area to volume ratio
$\frac{S(x)}{V(x)}$ decreases. In the case $x\in \mathbb{R}$, this
fact, equation (\ref{WaveFrontSpeed}) and Theorem 3.4 imply that the
wavefront of $u^{\epsilon,\delta}$ (for $\epsilon>0$ and $\delta>0$
small enough) slows down when the tube becomes thicker. A similar
result holds for layers.

\subsection{K-P-P Fronts in Random Media}
\vspace{0.2cm}

In this subsection we consider wave front propagation for the
solution of (\ref{InitialEquationWith_e}) for small
$\epsilon>0$, when $x\in \mathbb{R}$, the boundary $\partial D^{1}$
of $D^{1}$ is determined by stationary and ergodic random processes
on $\mathbb{R}$ and the nonlinear boundary term in
(\ref{InitialEquationWith_e}) (for $y=0$, i.e. $c(x,0,u)$) is of
K-P-P type. As we did in subsections 4.1, we will first
consider (Theorem 4.2.6) wavefront propagation for the solution of
(\ref{InitialEquationWithout_e}) and then with the aid of Theorem
3.4 we will consider (Theorem 4.2.7) wavefront propagation for the
solution of (\ref{InitialEquationWith_e}) for small enough
$\epsilon>0$. As we will see the cross sections $D_{x}$ of $D$
affect the speed of the wavefront through the surface to volume
ratio $\frac{S(x)}{V(x)}$.

In sections $7.4-7.6$ of \cite{K1} wave front propagation for
equations like (\ref{InitialEquationWithout_e}) is considered in the
case where there is no drift term and the randomness comes only from
the nonlinear part of the equation. Moreover in \cite{NX} the
authors considered the case of reaction-diffusion equations of type
(\ref{InitialEquationWithout_e}) with a random drift and homogeneous
in $x$ nonlinear term. In the case considered here, both the drift
and and the nonlinear term are random. In \cite{K1}, pp. 524-525,
the author remarks that one could use the procedure developed in
sections  $7.4-7.6$ of \cite{K1} to study wavefronts in
one-dimensional uniformly bounded random drift with random nonlinear
term. We will see that one can prove Theorem 4.2.6, which is
analogous to Theorem 7.6.1 in \cite{K1}, by following the proof of
Theorem 7.6.1 in \cite{K1}. We make use of the results in \cite{NX}
and of the fact that the operator of the equation
(\ref{InitialEquationWithout_e}) is self adjoint with respect to an appropriate inner product (it has the form
$\frac{1}{2V(x)}\frac{d}{dx}(V(x)\frac{d}{dx})$).  Actually the
latter simplifies the analysis significantly.  Instead of repeating
the proof of \cite{K1} here, we will only outline the differences.

Let us first list our assumptions. Consider a probability space
$(\hat{\Omega},\hat{\mathfrak{F}},\hat{P})$. We assume that the
random field $V(x,\hat{\omega})$ (namely the volume) is three times
continuously differentiable, i.e. $V\in
\mathcal{C}^{3}(\mathbb{R})$, with $\hat{P}$ probability one.
Suppose that $\Theta(x)=(\frac{d }{dx}(\log
V(x)),\frac{S(x)}{V(x)})$ is a random vector function on
$(\hat{\Omega},\hat{\mathfrak{F}},\hat{P})$ and that it is
measurable, stationary in $x$ and translation in $x$ generates an
ergodic transformation of the space $\hat{\Omega}$. Moreover the
function $\frac{d }{dx}(\log V(x))$ is assumed bounded, with zero
mean (i.e. $\hat{E}[\frac{d }{dx}(\log V(x))]=0$). We additionally
assume (for the purposes of Lemma 4.2.1 and 4.2.3) that there is a
set of nonzero $\hat{P}$ probability on which
\begin{equation}
\lim_{z\rightarrow
\infty}\int_{0}^{z}[V(x,\hat{\omega})]^{-1}dx=+\infty.\label{ConditionsForDrift}
\end{equation}
If condition (\ref{ConditionsForDrift}) holds on a set of nonzero
measure then, by the ergodicity assumption, it must hold with
$\hat{P}$ probability one.

As far as the non-linear term
$\overline{c}(x,u,\hat{\omega})u=\frac{S(x,\hat{\omega})}{V(x,\hat{\omega})}c(x,0,u)u$
is concerned, in addition to the stationarity and ergodicity
assumptions, we also make the following assumptions. For all $x\in
\mathbb{R}$, $c$ is of K.P.P type, i.e. $c(x,0,u)$ is positive for
$u<1$, negative for $u>1$, continuous in $u$ for $u\geq 0$ and
$c(x)=c(x,0,0)=\sup_{0 < u }c(x,0,u)$. Moreover with $\hat{P}$
probability one, the function $\overline{c}(x,u,\hat{\omega})u$
satisfies a Lipschitz condition of the form
$$|\overline{c}(x,u_{1},\hat{\omega})u_{1}-\overline{c}(x,u_{2},\hat{\omega})u_{2}|\leq
\frac{S(x,\hat{\omega})}{V(x,\hat{\omega})}\zeta(x)|u_{1}-u_{2}|,
\textrm{ for }x,u_{1},u_{2}\in \mathbb{R},$$
such that for all $t\geq 0$ and $x\in
\mathbb{R}$,
$$E_{x}\exp\{\int_{0}^{t}\frac{S(X_{s})}{V(X_{s})}\zeta(X_{s})\}<\infty, \hspace{0.2cm} \hat{P}\textrm{-a.s.},  $$

where $(X_{t},P_{x})$ is a diffusion process with random generator
$L=\frac{1}{2}\frac{d^{2}}{dx^{2}}+\frac{1}{2}\frac{d }{dx}(\log
V(x,\hat{\omega}))\frac{d}{dx}$.

The initial function $f(x)$ is assumed to be nonnegative, bounded from above and non-random.

\vspace{0.2cm}

Let now $\mu(z)$ be the function defined by the equality
\begin{equation}
\mu(z)=\hat{E}[\ln E_{1} \chi_{\tau_{0}<\infty}\exp\{\int_{0}^{\tau_{0}}[\overline{c}(X_{s})+z]ds\}], \hspace{0.2cm} z\in \mathbb{R},\label{MuFunction}
\end{equation}
where $\bar{c}(x)=\frac{1}{2}\frac{S(x)}{V(x)}c(x,0,0)$ and
$\tau_{0}$ is the first hitting time of the process $X_{t}$ to the
point $0$. For $\tau_{0}$ one has the following lemma:

\textbf{Lemma 4.2.1.} Condition (\ref{ConditionsForDrift}) and
$\hat{E}[\frac{d }{dx}(\log V(x,\hat{\omega}))]=0$ imply that
$P_{1}(\tau_{0}<\infty)=1$.

\vspace{0.2cm}

\textbf{Proof.}  It follows directly from the proof of Lemma 4.4 of
\cite{NX} if one notes that the drift term is $\frac{1}{2}\frac{d
}{dx}(\log V(x))$.

\begin{flushright}
$\square$
\end{flushright}

\vspace{0.2cm}

\textbf{Lemma 4.2.2.} Under the assumptions imposed above, function
$\mu(z)$ has the following properties:
\begin{enumerate}
\item{For all $z\in \mathbb{R}$, $\mu(z)=\lim_{t\rightarrow \infty}\frac{1}{t}\ln E_{t}\chi_{\tau_{0}<\infty}\exp\{\int_{0}^{\tau_{0}}[\overline{c}(X_{s})+z]ds\}$.}
\item{Function $\mu(z)$ is convex, lower
semicontinuous and monotonically non-decreasing in $z$. Moreover
$\mu(z)$ is continuously differentiable and the derivative $\mu'(z)$
is positive and monotonically increasing for $z<\bar{g}_{\mu}$,
where $\bar{g}_{\mu}$ is a non-positive number (which actually is
the discontinuity point of $\mu(z)$, as property (iii) below
shows).}
\item{$\mu(z)\leq 0$ for $z\leq \bar{g}_{\mu}$ and $\mu(z)=\infty$ for $z> \bar{g}_{\mu}$ where $\bar{g}_{\mu} \leq 0$.}
\end{enumerate}

\vspace{0.2cm}

\textbf{Proof.} Property (i) can be proven as Proposition 2.1 of
\cite{NX}. Property (ii) follows similarly as Theorem 7.5.1(ii) of
\cite{K1}. Property (iii) follows analogously to Theorem 7.5.1(iii)
of \cite{K1}. Here one uses the fact that the operator of
(\ref{InitialEquationWithout_e}) has the form
$\frac{1}{2V(x)}\frac{d}{dx}(V(x)\frac{d}{dx})$), i.e. it is self
adjoint.

\begin{flushright}
$\square$
\end{flushright}

We also observe that $\mu(z)\geq \mu_{o}(z)$ where
$\mu_{o}(z)=\hat{E}[\ln E_{1}(
\chi_{\tau_{0}<\infty}e^{z\tau_{0}})]$. As it has been proven in
Lemma 2.2 of \cite{NX}, function $\mu_{o}(z)$ has properties
(i)-(iii) of Lemma 4.2.2 as well (for $\bar{c}(x)=0$). In addition
the following lemma holds, which is a restatement of Proposition 4.1
of \cite{NX}.

\vspace{0.2cm}

\textbf{Lemma 4.2.3.} Condition (\ref{ConditionsForDrift}) and
$\hat{E}[\frac{d }{dx}(\log V(x,\hat{\omega}))]=0$ imply that the
discontinuity point of $\mu_{o}(z)$ is $\bar{g}_{\mu_{o}}=0$.

\vspace{0.2cm}

We will assume that $-\infty<\bar{g}_{\mu}<0$ (by Lemma 4.2.2(iii)
or Lemma 4.2.3 we already know that $\bar{g}_{\mu}\leq 0$) and we
define $I(y)=\sup_{z\leq \bar{g}_{\mu}}[yz-\mu(z)]$ for $y\in
\mathbb{R}$.

\vspace{0.2cm}

Lemma 4.2.2 and the fact that $\mu(z)\geq \mu_{o}(z)$ imply that the
arguments in the beginning of section $7.6$ of \cite{K1} carry out
here as well. Therefore we conclude that there is a unique
$\nu^{*}>0$ such that $I(\frac{1}{\nu^{*}})=0$ and
$\nu^{*}=\inf_{z\leq \bar{g}_{\mu}}\frac{z}{\mu(z)}$.

\vspace{0.2cm}

\textbf{Remark 4.2.4.} We would like to emphasize that the existence
and uniqueness of a positive $\nu^{*}$ follows mainly from
properties (i)-(iii) of $\mu(z)$ (Lemma 4.2.2). In particular
property (iii) holds because the operator of
(\ref{InitialEquationWithout_e}) is self adjoint.

\vspace{0.2cm}

Similarly as Theorem 7.6.1 in \cite{K1} was proven, one can prove
Theorem 4.2.6 below.

Note that by following the proof of Theorem 7.6.1 in \cite{K1}, one
needs to estimate certain probabilities for $\tau_{0}$ and $X_{t}$.
For this purpose we have the following lemma:
\vspace{0.2cm}

\textbf{Lemma 4.2.5.} Let $\delta$ be a positive number and
$U_{\delta}(0)=\lbrace x: |x|\leqslant \delta \rbrace$. Then
\begin{enumerate}
\item{ $\inf_{x\in U_{\delta}(0)} P_{x}\lbrace \tau_{0}\leqslant
1\rbrace>0$, $\hat{P}$-a.s.}
\item{$\inf_{x\in U_{\delta}(0), s\in(0,1]}
P_{x}\lbrace X_{s}\in U_{\delta}(0)\rbrace>0$, $\hat{P}$-a.s.}
\item{For $a>0$ and $\eta>\delta>0$ we have $$\inf_{x\in U_{\delta}(-a)}P_{x}\lbrace
\tau_{-\eta-a}> 1, X_{1}\in U_{\delta}(0)\rbrace>0,\hspace{0.1cm}
\hat{P}\textrm{-a.s.}$$ }
\end{enumerate}

\vspace{0.2cm}

\textbf{Proof.}  The proof of  all statements follows from the
corresponding statements for $W_{t}^{1}$ in place of $X_{t}$ (see
for example section 7.5 of \cite{K1}) and by the Girsanov's theorem
on the absolute continuous change of measures in the space of
trajectories.

\begin{flushright}
$\square$
\end{flushright}
Therefore we have the following Theorem:

\vspace{0.2cm}

\textbf{Theorem 4.2.6.} Let $x\in \mathbb{R}$ and $u(t,x)$ satisfy
equation (\ref{InitialEquationWithout_e}). Under our assumptions we
have:
\begin{enumerate}
\item{For all $\nu>\nu^{*}$,
\begin{displaymath}
\lim_{t\rightarrow\infty}\sup_{x\geq \nu t}u(t,x)=0, \hspace{0.2cm} \hat{P}-a.s.
\end{displaymath}
}
\item{Let us define  $\bar{c}_{h}(x)=\frac{1}{2}\frac{S(x)}{V(x)}\inf_{0< u < h}c(x,0,u)$ and assume that there is a constant $\kappa>0$ such that for any $0<h<1$ and $x\in \mathbb{R}$,
\begin{displaymath}
\kappa<\bar{c}_{h}(x), \hspace{0.2cm} \hat{P}-a.s.
\end{displaymath}
Then for all $\nu\in(0,\nu^{*})$,
\begin{displaymath}
\lim_{t\rightarrow\infty}\sup_{0\leq x\leq \nu t}u(t,x)=1, \hspace{0.2cm} \hat{P}-a.s.
\end{displaymath}
}
\end{enumerate}
Finally Theorem 3.4 and Theorem 4.2.6 imply:

\vspace{0.2cm}

\textbf{Theorem 4.2.7.} Let $(x,y)\in \mathbb{R}\times
\mathbb{R}^{m}$ and $u^{\epsilon}(t,x,y)$ satisfy equation
(\ref{InitialEquationWith_e}). Under our assumptions we have:
\begin{enumerate}
\item{For all $\nu>\nu^{*}$,
\begin{displaymath}
\lim_{t\rightarrow\infty}\sup_{x\geq \nu t}\lim_{\epsilon\rightarrow 0}u^{\epsilon}(t,x,y)=0, \hspace{0.2cm} \hat{P}-a.s.
\end{displaymath}
}
\item{Let us define  $\bar{c}_{h}(x)=\frac{1}{2}\frac{S(x)}{V(x)}\inf_{0< u < h}c(x,0,u)$ and assume that there is a constant $\kappa>0$ such that for any $0<h<1$ and $x\in \mathbb{R}$,
\begin{displaymath}
\kappa<\bar{c}_{h}(x), \hspace{0.2cm} \hat{P}-a.s.
\end{displaymath}
Then for all $\nu\in(0,\nu^{*})$,
\begin{displaymath}
\lim_{t\rightarrow\infty}\sup_{0\leq x\leq \nu t}\lim_{\epsilon\rightarrow 0}u^{\epsilon}(t,x,y)=1, \hspace{0.2cm} \hat{P}-a.s.
\end{displaymath}
}
\end{enumerate}

\vspace{0.2cm}

\textbf{Remark 4.2.8.} Theorem 4.2.6 was proven in (\cite{K1}) with
the assumption in part (ii) replaced by the assumption that for any
$0<h<1$ and $\nu \in \mathbb{R}$,
\begin{equation}
\limsup_{t\rightarrow\infty}\frac{1}{t}\ln E_{\nu
t}\exp\{-\int_{0}^{t}\bar{c}_{h}(X_{s})ds\}<0, \hspace{0.2cm}
\hat{P}-a.s,\label{Assumption1}
\end{equation}
which is however difficult to verify. Obviously the assumption made
in part (ii) of Theorems 4.2.6 and 4.2.7 implies
(\ref{Assumption1}).

\section{Acknowledgements}
M. Freidlin was partially supported by the NSF. The authors would like to thank Prof. K. Trivisa and Prof. M. Grillakis for helpfull discussions.


\end{document}